\documentclass{article}

\input{diagrams}

\usepackage{amstex}
\usepackage{latexsym}

\newcommand{\mcm}[3]{\newcommand{#1}[#2]{{\ensuremath{#3}}}}

\mcm{\mc}{1}{\mathcal{#1}}
\mcm{\mr}{1}{\mathrm{#1}}
\mcm{\mi}{1}{\mathit{#1}}
\mcm{\mb}{1}{\mathbf{#1}}
\mcm{\cat}{1}{\mc{#1}}
\mcm{\scat}{1}{\Bbb{#1}}
\mcm{\fcat}{1}{\mb{#1}}  	
\mcm{\twid}{1}{\widetilde{#1}}

\mcm{\pr}{2}{(#1,#2)}
\mcm{\triple}{3}{(#1,#2,#3)}
\mcm{\range}{2}{#1,\,\ldots\,,#2}
\mcm{\tuplebts}{1}{\langle #1 \rangle}
\mcm{\tuple}{3}{\tuplebts{\range{#1,#2}{#3}}}
\mcm{\bftuple}{2}{\tuplebts{\range{#1}{#2}}}
\mcm{\rttuple}{1}{\tuplebts{\,\ldots\,,#1}}
\mcm{\elt}{0}{\in}
\mcm{\such}{0}{\:|\:}
\mcm{\nat}{0}{\Bbb{N}}	
\mcm{\iso}{0}{\cong}
\mcm{\eqv}{0}{\simeq}
\mcm{\id}{0}{\mi{id}}
\mcm{\of}{0}{\circ}
\mcm{\blank}{0}{(\:\:)}
\mcm{\dashbk}{0}{\mbox{\underline{\hspace{1em}}}}
\mcm{\op}{0}{\mr{op}}
\mcm{\ftrcat}{2}{[#1,#2]}
\mcm{\Sets}{0}{\fcat{Sets}}
\newcommand{\pf}{\textbf{Proof}}
\newcommand{\done}{\hfill\ensuremath{\Box}}

\mcm{\go}{0}{\rTo}
\mcm{\og}{0}{\lTo}
\mcm{\goby}{1}{\rTo^{#1}}
\mcm{\ogby}{1}{\lTo^{#1}}
\mcm{\goesto}{0}{\,\longmapsto\,}
\mcm{\slob}{3}{(#1 \goby{#2} #3)}
\mcm{\bktdslob}{3}{\slob{#1}{#2}{#3}}

\mcm{\ess}{0}{\cat{S}}
\mcm{\blob}{0}{\scriptscriptstyle{\bullet}}
\mcm{\ust}{0}{{}^{*}}
\mcm{\ubl}{0}{{}^{\blob}}
\mcm{\stbk}{0}{\blank^{*}}
\mcm{\blbk}{0}{\blank^{\blob}}
\mcm{\Mnd}{0}{\triple{\stbk}{\eta}{\mu}}
\mcm{\Cartpr}{0}{\pr{\cat{S}}{\ust}}
\mcm{\dom}{0}{\mr{dom}}
\mcm{\cod}{0}{\mr{cod}}
\mcm{\spn}{3}{#2 \og #1 \go #3}
\mcm{\spaan}{5}{#2 \ogby{#4} #1 \goby{#5} #3}
\mcm{\gph}{2}{\spn{#1}{#2^{*}}{#2}}
\mcm{\graph}{4}{\spaan{#1}{#2^{*}}{#2}{#3}{#4}}
\mcm{\ids}{0}{\mi{ids}}
\mcm{\comp}{0}{\mi{comp}}
\mcm{\Mon}{0}{\mb{Mon}}
\mcm{\Slice}{0}{\cat{S}/S}
\mcm{\unit}{0}{\mi{unit}}
\mcm{\mult}{0}{\mi{mult}}
\mcm{\Imnd}{0}{\triple{\blbk}{\unit}{\mult}}
\mcm{\Icartpr}{0}{\pr{\Slice}{\ubl}}
\mcm{\Alg}{0}{\mb{Alg}}
\mcm{\mtihom}{4}{#1(\range{#2}{#3};#4)} 
\mcm{\mtiendo}{2}{#1(#2;#2)}

\newcommand{\pullshape}
	{\setlength{\unitlength}{1em}
	\begin{picture}(2,5)(-1,-5)
	\put(0,-5){\line(1,1){1}}
	\put(0,-5){\line(-1,1){1}}
	\end{picture}}
\newcommand{\Spbk}{\overprint{\raisebox{-2.5em}{\pullshape}}}
\newenvironment{slopeydiag}
	{\begin{diagram}[size=2em]}
	{\end{diagram}}
\mcm{\diso}{0}{\sim}     
\mcm{\vdiso}{0}{\wr}
\newarrow{Get}....>
\newenvironment{triangdiag}
	{\begin{diagram}[width=1em,height=1.5em]}
	{\end{diagram}}

\mcm{\node}{0}{\bullet}
\mcm{\enode}{0}{\circ}
\mcm{\nl}{1}{\stackrel{#1}{\node}}
\newcommand{\nodelt}[1]{\raisebox{1ex}{\ensuremath{#1}}}
\newenvironment{tree}
	{\begin{diagram}[height=1em,width=.75em,abut,noPS]}	
	{\end{diagram}}
\newcommand{\lt}[1]{\ldLine(#1,2)}
\newcommand{\rt}[1]{\rdLine(#1,2)}
\newcommand{\dn}{\dLine}
\mcm{\nll}{1}{\stackrel{\node}{#1}}

\mcm{\Hom}{0}{\mr{Hom}}
\mcm{\homset}{3}{#1(#2,#3)}
\mcm{\End}{0}{\mr{End}}
\mcm{\Sym}{0}{\mr{Sym}}
\mcm{\vslob}{3}
	{\left.
	\begin{diagram}[height=1.5em]
	#1		\\
	\dTo>{#2}	\\
	#3		\\
	\end{diagram}
	\right.}
\mcm{\bktdvslob}{3}{\left( \vslob{#1}{#2}{#3} \right)}
\mcm{\emptybk}{0}{\:\:}
\newenvironment{ntdiag}
	{\begin{diagram}[size=1.5em,noPS]}
	{\end{diagram}}
\mcm{\Cat}{0}{\fcat{Cat}}
\mcm{\Multicat}{0}{\fcat{Multicat}}
\mcm{\Multicatof}{1}{#1\!\!-\!\!\Multicat}
\mcm{\Graph}{0}{\fcat{Graph}}
\mcm{\Graphof}{1}{#1\!\!-\!\!\Graph}
\mcm{\lid}{0}{\star}
\mcm{\ulid}{0}{{}^\lid}
\mcm{\lidbk}{0}{\blank^{\lid}}
\newcommand{\url}[1]{\mbox{\tt #1}}
\mcm{\without}{0}{\setminus}
\mcm{\sub}{0}{\subseteq}
\mcm{\ush}{0}{{}^{\sharp}}
\newarrow{Edge}---->
\newcommand{\cnr}{}
\newenvironment{opetope}
	{\begin{diagram}[size=1em,abut,tight,noPS]}
	{\end{diagram}}

\newcommand{\EmptyOne}
{\begin{tree}
\enode	\\
\dn	\\
\node	\\
\end{tree}}

\newcommand{\Oak}[5]
{\begin{tree}
 & & & &\enode& & & & & & & & \\
 & & & &\dn& & & & & & & & \\
\nl{#1}& &\nl{#2}& &\node& & & &\nl{#4}& & & &\nl{#5}\\
 &\rt{2}&\dn&\lt{2}& & & & & &\rt{2}& &\lt{2}& \\
 & &\node& & & &\nl{#3}& & & &\node& & \\
 & & &\rt{4}& & &\dn& & &\lt{4}& & & \\
 & & & & & &\node& & & & & & \\
\end{tree}}

\newcommand{\Pear}[2]
{\begin{tree}
	&	&	&	&\enode	&&	&	&	\\
	&	&	&	&\dn	&&	&	&	\\
\nl{#1}&	&	&	&\node	&&	&	&	\\
	&\rt{2}	&	&\lt{2}	&	&&	&	&	\\
	&	&\node	&	&	&&	&	&\nl{#2}\\
	&	&	&\rt{4}	&	&&	&\lt{2}	&	\\
	&	&	&	&	&&\node	&	&	\\
\end{tree}}

\newcommand{\Orange}[2]
{\begin{tree}
	&	&	&	&\nl{#2}\\
	&	&	&	&\dn	\\
\nl{#1}&	&	&	&\node	\\
	&\rt{2}	&	&\lt{2}	&	\\
	&	&\node	&	&	\\
\end{tree}}

\newcommand{\Apple}[2]
{\begin{tree}
	&	&\enode	&	&	\\
	&	&\dn	&	&	\\
\nl{#1}	&	&\node	&	&\nl{#2}\\
	&\rt{2}	&\dn	&\lt{2}	&	\\
	&	&\node	&	&	\\
	&	&\dn	&	&	\\
	&	&\node	&	&	\\
\end{tree}}

\newcommand{\MixedFruit}[4]
{\begin{tree}
 & & & &\nl{#2}& & & & & & &\enode& & \\
 & & & &\dn    & & & & & & &\dn   & & \\
\nl{#1}& & & &\node& & &\enode& &\nl{#3}& &\node& &\nl{#4}\\
 &\rt{2}& &\lt{2}& & & &\dn& & &\rt{2}&\dn&\lt{2}& \\
 & &\node& & & & &\node& & & &\node& & \\
 & & &\rt{3}& & &\lt{2}& & & & &\dn& & \\
 & & & & &\node& & & & & &\node& & \\
 & & & & & &\rt{3}& & & &\lt{3}& & & \\
 & & & & & & & &\node& & & & & \\
\end{tree}}

\newlength{\volt}
\newcommand{\transistor}[5]
{\setlength{\unitlength}{1\volt}
\begin{picture}(18,12)(-5,-6)
\put(0,6){\line(0,-1){12}}
\put(0,-6){\line(3,2){9}}
\put(0,6){\line(3,-2){9}}
\put(9,0){\line(1,0){2}}
\put(-2,4){\line(1,0){2}}
\put(-2,2){\line(1,0){2}}	
\put(-2,-4){\line(1,0){2}}
\put(12,-0.5){\ensuremath{#5}}
\put(-5,3.5){\ensuremath{#2}}
\put(-5,1.5){\ensuremath{#3}}		
\put(-5,-4.5){\ensuremath{#4}}
\thicklines
\put(-1.5,0){\line(1,0){.1}}	
\put(-1.5,-1){\line(1,0){.1}}		
\put(-1.5,-2){\line(1,0){.1}}		
\thinlines
\put(2,-0.5){\ensuremath{#1}}
\end{picture}}
\newcommand{\ctransistor}[5]
	{\raisebox{-6\volt}{\transistor{#1}{#2}{#3}{#4}{#5}}}

\newcommand{\bftransistor}[4]
{\setlength{\unitlength}{1\volt}
\begin{picture}(18,12)(-5,-6)
\put(0,6){\line(0,-1){12}}
\put(0,-6){\line(3,2){9}}
\put(0,6){\line(3,-2){9}}
\put(9,0){\line(1,0){2}}
\put(-2,4){\line(1,0){2}}
\put(-2,-4){\line(1,0){2}}
\put(12,-0.5){\ensuremath{#4}}
\put(-5,3.5){\ensuremath{#2}}
\put(-5,-4.5){\ensuremath{#3}}
\thicklines
\put(-1.5,1){\line(1,0){.1}}		
\put(-1.5,0){\line(1,0){.1}}		
\put(-1.5,-1){\line(1,0){.1}}
\thinlines		
\put(2,-0.5){\ensuremath{#1}}
\end{picture}}
\newcommand{\comptrans}[4]
{\setlength{\volt}{.5ex}
\setlength{\unitlength}{1\volt}
\begin{picture}(80,72)(0,-36)
\put(0,23){\bftransistor{#1}{}{}{}}
\put(0,6){\bftransistor{#2}{}{}{}}
\put(0,-35){\bftransistor{#3}{}{}{}}
\put(63,-6){\transistor{#4}{}{}{}}
\put(17,29){\line(2,-1){50}}
\put(17,12){\line(5,-1){50}}
\put(17,-29){\line(2,1){50}}
\thicklines
\put(24,-5){\line(1,0){.4}}		
\put(24,-9){\line(1,0){.4}}		
\put(24,-13){\line(1,0){.4}}
\thinlines		
\setlength{\volt}{1ex}
\end{picture}}

\newcommand{\threetransistor}[5]
{\setlength{\unitlength}{1\volt}
\begin{picture}(18,12)(-5,-6)
\put(0,6){\line(0,-1){12}}
\put(0,-6){\line(3,2){9}}
\put(0,6){\line(3,-2){9}}
\put(9,0){\line(1,0){2}}
\put(-2,4){\line(1,0){2}}
\put(-2,0){\line(1,0){2}}	
\put(-2,-4){\line(1,0){2}}
\put(12,-0.5){\ensuremath{#5}}
\put(-5,3.5){\ensuremath{#2}}
\put(-5,-0.5){\ensuremath{#3}}		
\put(-5,-4.5){\ensuremath{#4}}
\put(2,-0.5){\ensuremath{#1}}
\end{picture}}
\newcommand{\twotransistor}[4]
{\setlength{\unitlength}{1\volt}
\begin{picture}(18,12)(-5,-6)
\put(0,6){\line(0,-1){12}}
\put(0,-6){\line(3,2){9}}
\put(0,6){\line(3,-2){9}}
\put(9,0){\line(1,0){2}}
\put(-2,4){\line(1,0){2}}
\put(-2,-4){\line(1,0){2}}
\put(12,-0.5){\ensuremath{#4}}
\put(-5,3.5){\ensuremath{#2}}
\put(-5,-4.5){\ensuremath{#3}}
\put(2,-0.5){\ensuremath{#1}}
\end{picture}}
\newcommand{\notransistor}[2]
{\setlength{\unitlength}{1\volt}
\begin{picture}(18,12)(-5,-6)
\put(0,6){\line(0,-1){12}}
\put(0,-6){\line(3,2){9}}
\put(0,6){\line(3,-2){9}}
\put(9,0){\line(1,0){2}}
\put(12,-0.5){\ensuremath{#2}}
\put(2,-0.5){\ensuremath{#1}}
\end{picture}}
\newcommand{\freemoncatpic}
{\setlength{\volt}{.6ex}
\setlength{\unitlength}{1\volt}
\begin{picture}(18,48)
\put(0,3){\twotransistor{a_3}{s_4}{s_5}{s'_3}}
\put(0,18){\notransistor{a_2}{s'_2}}
\put(0,33){\threetransistor{a_1}{s_1}{s_2}{s_3}{s'_1}}
\put(5.7,0){\framebox(10.3,48){}}
\end{picture}}
\newcommand{\dotty}{\raisebox{-.4ex}{\ensuremath{\!\scriptstyle\bullet}}}
\newcommand{\discfibpic}
{\setlength{\unitlength}{1em}		
\begin{picture}(16,9)(-6,0)
\put(4,2){\oval(8,4)}
\put(4,8){\dotty}
\put(2,2){\dotty}
\put(4,1){\dotty}
\put(6,1){\dotty}
\put(5,3){\dotty}
\put(4,8){\vector(-1,-3){2}}
\put(2,2){\vector(2,-1){2}}
\put(4,8){\vector(0,-1){7}}
\put(2,2){\vector(3,1){3}}
\put(4,1){\vector(1,0){2}}
\put(4.2,8){0}
\put(4.2,5.5){\ensuremath{f\of a}}
\put(2.5,5.5){\ensuremath{a}}
\put(1.2,1.7){\ensuremath{s}}
\put(3.6,.2){\ensuremath{s'}}
\put(2.1,.9){\ensuremath{f}}
\put(8.5,1.5){\ensuremath{S}}
\end{picture}}
\newcommand{\cdiscfibpic}{\raisebox{-4.5em}{\discfibpic}}
\newtheorem{thm}{Theorem}[section]
\newtheorem{defn}[thm]{Definition}

\newtheorem{lemma}[thm]{Lemma}
\newtheorem{constn}[thm]{Construction}

\newtheorem{lotsofremarks}[thm]{Remarks}
\newenvironment{remarks}[1]
	{\begin{lotsofremarks}\label{#1}\end{lotsofremarks}\begin{enumerate}}
	{\end{enumerate}}

\newtheorem{lotsofegs}[thm]{Examples}
\newenvironment{eg}[1]
	{\begin{lotsofegs}\label{#1}\end{lotsofegs}\begin{enumerate}}
	{\end{enumerate}}
\newtheorem{discussion}[thm]{Discussion}

\newtheorem{llemma}{Lemma}[subsection]

\title{General Operads and Multicategories}
\author{Tom Leinster\\ \\
	\normalsize{Department of Pure Mathematics, University of
	Cambridge}\\ 
	\normalsize{Email: leinster@@dpmms.cam.ac.uk}\\
	\normalsize{Web: http://www.dpmms.cam.ac.uk/$\sim$leinster}}
\date{12 October, 1997\\
	\small Minor amendments Dec 97, Oct 98}

\begin{document}


\begin{titlepage}

\maketitle
\thispagestyle{empty}

\begin{abstract}
Notions of `operad' and `multicategory' abound. This work provides a single
framework in which many of these various notions can be expressed.
Explicitly: given a monad \ust\ on a category \ess, we define the term
\emph{\Cartpr-multicategory}, subject to certain conditions on \ess\ and
\ust.  Different choices of \ess\ and \ust\ give some of the existing
notions. We then describe the \emph{algebras} for an \Cartpr-multicategory
and, finally, present a tentative selection of further developments. Our
approach makes possible concise descriptions of Baez and Dolan's opetopes and
Batanin's operads; both of these are included.
\end{abstract}

\vspace{2mm}

\begin{center} 
\bfseries Contents
\end{center} 

\begin{center} 
\begin{tabular}{rlcr}
	&Introduction			&&\pageref{sec:intro}		\\
1	&Cartesian Monads		&&\pageref{sec:cart-mnds}	\\
2	&Multicategories		&&\pageref{sec:multicats}	\\
3	&Algebras			&&\pageref{sec:algs}		\\
4	&Further Developments		&&\pageref{sec:fds}		\\
	&References			&&\pageref{sec:biblio}
\end{tabular}
\end{center}

\end{titlepage}



\section*{Introduction}	\label{sec:intro}

Operads, multicategories and the like have appeared in many different guises,
especially recently. Investigators of $n$-categories (Baez-Dolan, Batanin,
Hermida-Makkai-Power) have had cause to resurrect and generalize in various
ways the original May definition of c.1970. The work of Soibelman and
of Borcherds, related to topological quantum field theory and vertex
algebras, calls upon still different ideas of multicategory.

It is the aim of this article to unify some, if not all, of these approaches.
Where it works, it provides a single formalism in which these various notions
can be expressed simply and perhaps compared. It does succeed in capturing
the `traditional' definitions of (non-symmetric) operad and multicategory,
the Batanin operads, and at least some of the flavour of Baez-Dolan,
Soibelman and Borcherds. (An optimist might even imagine that it would
facilitate the comparison of different definitions of weak $n$-category.) A
notable failure is that the symmetric group actions often included in
definitions of `operad' are not easily expressed in our system, for the time
being at least.

The central idea is simple. In a plain category, an arrow $a$ is written
\[
s' \goby{a} s,
\]
where $s'$ and $s$ are elements of the set $S$ of objects. In an `ordinary'
multicategory, an arrow $a$ is written
\[
\range{s_1}{s_n} \goby{a} s,
\]
where $s$ is an element of the set $S$ of objects, \bftuple{s_1}{s_n} is an
element of $S^*$, and \ust\ is the free-monoid monad (`word monad') on \Sets.
Thus the graph structure of the multicategory is a diagram
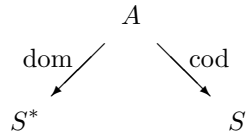
\begin{slopeydiag}
	&		&A	&		&	\\
	&\ldTo<{\dom}	&	&\rdTo>{\cod}	&	\\
S^{*}	&		&	&		&S	\\
\end{slopeydiag}
in \Sets, where $A$ is the set of arrows. Now, just as a (small) category can
be described as a diagram 
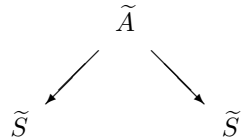
\begin{slopeydiag}
		&	&\twid{A}	&	&		\\
		&\ldTo	&		&\rdTo	&		\\
\twid{S}	&	&		&	&\twid{S}	\\
\end{slopeydiag}
in \Sets\ together with identity and composition functions
\[
\twid{S}\go\twid{A} \mr{,\ \ } \twid{A}\times_{\twid{S}}\twid{A}\go\twid{A}
\]
satisfying some axioms, so we may describe the multicategory structure on
$\left.
\begin{slopeydiag}
	&	&A	&	&	\\
	&\ldTo	&	&\rdTo	&	\\
S^{*}	&	&	&	&S	\\
\end{slopeydiag}
\right.$
by manipulation of certain diagrams in \Sets. In general, we take a category
\ess\ and a monad \ust\ on \ess, and, subject to certain conditions, define
`\Cartpr-multicategory'. Thus a category is a
\pr{\Sets}{\id}-multicategory, and an ordinary multicategory is a
\pr{\Sets}{\mr{free\ monoid}}-multicategory.

Section~\ref{sec:cart-mnds} describes the simple conditions needed on \ess\
and \ust\ in order that everything that follows will work. Many examples are
given. Section~\ref{sec:multicats} explains what an \Cartpr-multicategory is
and how the examples relate to existing notions of multicategory. In
particular, a concise definition of Batanin operads is given. Most of
these existing notions carry with them the concept of an \emph{algebra} for
an operad/multicategory; section~\ref{sec:algs} defines algebras in our
general setting. Finally, this being work in motion, section~\ref{sec:fds} is
a collection of loose topics and possible further developments, of assorted
merit. Included is a compact construction of Baez and Dolan's opetopes.

\subsection*{Related Work}
Since the original posting of this document, various pieces of related work
have been pointed out to me. The basic construction of
\Cartpr-multicategories was carried out by Burroni (\cite{Bur}) in 1971,
under the (better?) name of $T$-categories, where $T=\ust$ is the monad
concerned. He develops the theory extensively, although in what direction I
do not yet know. This reference was passed on to me by Claudio Hermida, who
also, independently, made this definition (and more). Notes from talks he has
given on the subject are available (\cite{Her}). The present work is also
connected to Kelly's theory of clubs, for which see \cite{Kel1} and
\cite{Kel2}.

\subsection*{Acknowledgements}
This work was supported by a PhD grant from EPSRC. It was prepared in
\LaTeX\ using Paul Taylor's commutative diagrams package. My thanks go to
Martin Hyland for encouraging me to pursue this line of enquiry. I am also
grateful to Peter Johnstone, Craig Snydal and Robin Cockett, for help and
useful discussions.


\clearpage



\section{Cartesian Monads}	\label{sec:cart-mnds}

In this section we introduce the conditions required on a monad \Mnd\ on a
category \ess, in order that we may (in Section~\ref{sec:multicats}) define
the notion of an \Cartpr-multicategory. Like most conditions that follow, the
demand is that certain things are or preserve pullbacks.

\begin{defn}
\label{defn-cart}
A monad \Mnd\ on a category \cat{S} will be called \emph{cartesian} if
\begin{enumerate}
\item $\eta$ and $\mu$ are cartesian natural transformations, i.e.
for any $X \goby{f} Y$ in \cat{S} the naturality squares
\[\left.
\begin{diagram}
X	&\rTo^{\eta_{X}}&X^{*}		\\
\dTo<{f}&		&\dTo>{f^{*}}	\\
Y	&\rTo^{\eta_{Y}}&Y^{*}		\\
\end{diagram}
\right.
\mi{and}
\left.
\begin{diagram}
X^{**}		&\rTo^{\mu_{X}}	&X^{*}		\\
\dTo<{f^{**}}	&		&\dTo>{f^{*}}	\\
Y^{**}		&\rTo^{\mu_{Y}}	&Y^{*}		\\
\end{diagram}
\right.\]
are pullbacks, and \label{nts-cart}
\item \stbk\ preserves pullbacks. \label{pb-pres}
\end{enumerate}
\end{defn}

(Conditions~\ref{nts-cart} and~\ref{pb-pres} ensure that not just $\eta$ and
$\mu$, but \emph{all} natural transformations $\blank^{*n} \go \blank^{*m}$
arising from the monad, are cartesian. For instance, $\mu*\eta^{*}:
\blank^{***} \go \blank^{***}$ is cartesian.)

Checking that a particular monad is cartesian can be eased slightly by
employing:

\begin{lemma}
\label{conv-crit}
Let \cat{S} be a category with a terminal object 1. Then a monad \Mnd\
on \cat{S} satisfies condition~\ref{nts-cart} of
Definition~\ref{defn-cart} iff for any object $Z$ of \cat{S}, the squares
\[\left.
\begin{diagram}
Z	&\rTo^{\eta_{Z}}&Z^{*}		\\
\dTo<{!}&		&\dTo>{!^{*}}	\\
1	&\rTo^{\eta_{1}}&1^{*}		\\
\end{diagram}
\right.
\mi{and}
\left.
\begin{diagram}
Z^{**}		&\rTo^{\mu_{Z}}	&Z^{*}		\\
\dTo<{!^{**}}	&		&\dTo>{!^{*}}	\\
1^{**}		&\rTo^{\mu_{1}}	&1^{*}		\\
\end{diagram}
\right.\]
are pullbacks, where $!$ is the unique map $Z \go 1$.
\end{lemma}

\pf\ Just use the fact that if in a commutative diagram
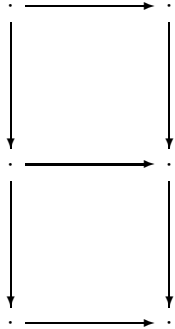
\begin{diagram}
\cdot	&\rTo	&\cdot	\\
\dTo	&	&\dTo	\\
\cdot	&\rTo	&\cdot	\\
\dTo	&	&\dTo	\\
\cdot	&\rTo	&\cdot	\\
\end{diagram}
the outer rectangle and lower square are pullbacks, then so too is the
upper square.~\done

The condition that \ess\ must satisfy is:

\begin{defn}
A category is called \emph{cartesian} if it has all finite limits.
\end{defn}

\begin{eg}{eg:cart-mnds}

\item The identity monad on any category is clearly cartesian.

\item \label{free-monoid-is-cart}
Let $\ess = \Sets$ and let \ust\ be the monoid monad, i.e. the
monad arising from the adjunction
\begin{diagram}
\fcat{Monoids}	&\pile{\rTo \\ \ \ \top \\ \lTo}	&\Sets.
\end{diagram}
Certainly \ess\ is cartesian. We show that \ust, too, is cartesian, using
Lemma~\ref{conv-crit}.

\emph{Unit:} Observe that $1^* = \nat$, that $1\goby{\eta_1}\nat$ has image
$\{1\}$, and that $X^*\goby{!^*}\nat$ sends a word
\tuple{x_1}{x_2}{x_n}\ to $n\elt\nat$. The square
\begin{diagram}
X	&\rTo^{\eta_X}	&X^*		\\
\dTo<!	&		&\dTo>{!^*}	\\
1	&\rTo^{\eta_1}	&1^*		\\
\end{diagram}
is a pullback, as $X \iso \{\bftuple{x_1}{x_n}\elt X^{*} \such n=1\}$.

\emph{Multiplication:} The map $\nat^{*}\goby{\mu_1}\nat$ is
$\bftuple{n_1}{n_k} \goesto n_{1}+ \cdots +n_{k},$ and 
$X^{**}\goby{!^{**}}\nat^{*}$ is 
$\bftuple{\bftuple{x^{1}_{1}}{x^{1}_{n_1}}}{\bftuple{x^{k}_{1}}{x^{k}_{n_k}}}
\goesto \bftuple{n_1}{n_k}.$ The square
\begin{diagram}
X^{**}	&\rTo^{\mu_X}	&X^*		\\
\dTo<{!^{**}}	&	&\dTo>{!^*}	\\
\nat^*	&\rTo^{\mu_1}	&\nat		\\
\end{diagram}
is a pullback: for
\begin{eqnarray*}
X^{*} \times_{\nat} \nat^*	
	&\iso	&\{\pr{\bftuple{x_1}{x_m}}{\bftuple{n_1}{n_k}} \such
			m=n_{1}+ \cdots +n_{k}\}	\\
	&\iso	&\{\tuple{\bftuple{x_1}{x_{n_1}}}
		         {\bftuple{x_{n_{1}+1}}{x_{n_{1}+n_{2}}}}
	      	         {\rttuple{x_{n_{1} + \cdots + n_{k}}}}\}	\\
	&\iso	&X^{**}.
\end{eqnarray*}

\emph{Pullback Preservation:} Let
\begin{diagram}
P	&\rTo^{}	&Y	\\
\dTo<{}	&		&\dTo>{g}	\\
X	&\rTo_{f}	&Z	\\
\end{diagram}
be a pullback square: then
\begin{eqnarray*}
X^{*} \times_{Z^*} Y^{*}	
	&\iso	&\{\pr{\bftuple{x_1}{x_n}}{\bftuple{y_1}{y_m}} \such
			\bftuple{fx_1}{fx_n}=\bftuple{gy_1}{gy_m}\}	\\
	&\iso	&\{\bftuple{\pr{x_1}{y_1}}{\pr{x_n}{y_n}} \such
			fx_{i} = gy_{i} \mr{\ for\ each\ } i\}	\\
	&\iso	&P^{*}.
\end{eqnarray*}

\item \label{eg:comm-not-cart}
A non-example. Let $\ess=\Sets$ and let \Mnd\ be the free commutative monoid
monad. This is not cartesian: e.g.\ the naturality square for $\mu$ at $2\go
1$ is not a pullback. 

\item \label{eg:sr-theories}
Let $\ess = \Sets$. Any finitary algebraic theory gives a monad
on \ess; which are cartesian? Without answering this question
completely, we indicate a certain class of theories which do give
cartesian monads. An equation (made up of variables and finitary
operators) is said to be \emph{strongly regular} if the same variables
appear in the same order, without repetition, on each side. Thus
\[\begin{array}{ccc}
(x.y).z=x.(y.z) & \mr{and} & (x\uparrow y)\uparrow z = x\uparrow(y.z),\\
\end{array}\]
but not
\[\begin{array}{cccc}
x+(y+(-y))=x, & x.y=y.x & \mr{or} & (x.x).y=x.(x.y),\\
\end{array}\]
qualify. A theory is called \emph{strongly regular} if it can be presented by
operators and strongly regular equations. It may be apparent from
Example~\ref{free-monoid-is-cart} that the only property we used of the
theory of monoids was its strong regularity, and that in a similar way we
could prove that any strongly regular theory yields a cartesian monad.

This last result, and the notion of strong regularity, are due to
Carboni and Johnstone. They show in \cite[Proposition 3.2 via Theorem 2.6]{CJ}
that a theory is strongly regular iff $\eta$ and $\mu$ are
cartesian natural transformations and $\stbk$ preserves wide
pullbacks. A \emph{wide pullback} is by definition a limit of shape
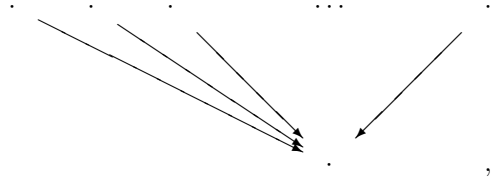
\begin{diagram}
\cdot	&\cdot	&\cdot		&	&\cdots	&		&\cdot	\\
	&\rdTo(4,2)&\rdTo(3,2)	&\rdTo(2,2)&	&\ldTo(2,2)	&	\\
	&	&		&	&\cdot	&		&,	\\
\end{diagram}
where the top row is a set of any size (perhaps infinite). When the set is of
size 2 this is an ordinary pullback, so the monad from a strongly regular
theory is indeed cartesian. (Examples~\ref{eg:exceptions-mnd},
\ref{eg:unary-mnd}, and \ref{eg:tree-mnd} can also be found in \cite{CJ}.)

\item \label{eg:exceptions-mnd}
Let $\ess = \Sets$, let $E$ be a fixed set, and let + denote binary
coproduct: then the endofunctor $\dashbk +E$ on \ess\ has a natural monad
structure. This monad is cartesian, corresponding to the algebraic theory
consisting only of one constant for each member of $E$.

\item \label{eg:unary-mnd}
Let $\ess=\Sets$ and let $M$ be a monoid: then the endofunctor
$M\times\dashbk$ on \ess\ has a natural monad structure. This monad is
cartesian, corresponding to an algebraic theory consisting only of unary
operations. 

\item \label{eg:tree-mnd}
Let $\ess = \Sets$, and consider the finitary algebraic theory
on $\ess$ generated by one $n$-ary operation for each $n\elt\nat$, and
no equations. This theory is strongly regular, so the induced monad
\Mnd\ on $\ess$ is cartesian.

If $X$ is any set then $X^*$ can be described inductively by:
\begin{itemize}
	\item if $x\elt X$ then $x\elt X^*$
	\item if $\range{t_1}{t_n} \elt X^*$ then
$\bftuple{t_1}{t_n}\elt X^*$.
\end{itemize}
We can draw any element of $X^*$ as a tree with leaves labelled by
elements of $X$:
\begin{itemize}
	\item $x\elt X$ is drawn as \nl{x}
	\item if \range{t_1}{t_n} are drawn as \range{T_1}{T_n} then
\bftuple{t_1}{t_n} is drawn as
$\left.
\begin{tree}
\nodelt{T_1}&	&\nodelt{T_2}&	&	&\cdots	& & &\nodelt{T_n}\\
	&\rt{4}	&	&\rt{2}	&	&	& &\lt{4} &	\\
	&	&	&	&\node	&	& & &		\\
\end{tree}
\right.$
, or if $n=0$, as $\left.\EmptyOne\right.$.
\end{itemize}
Thus the element $\tuplebts{\tuplebts{x_{1},x_{2},\tuplebts{}}, x_{3},
\tuplebts{x_{4},x_{5}}}$ of $X^*$ is drawn as
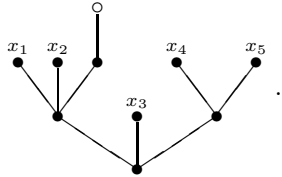
\[\left.\Oak{x_1}{x_2}{x_3}{x_4}{x_5}\right. .\]

The unit $X\go X^*$ is $x\goesto\nl{x}$, and multiplication $X^{**}
\go X^{*}$ takes an $X^*$-labelled tree (e.g.\ 
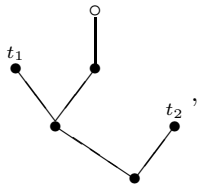
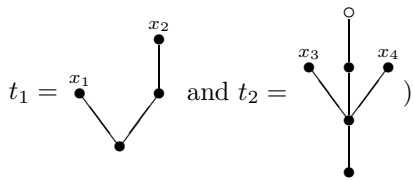
\[\left.\Pear{t_1}{t_2}\right.,\]
with
\[t_{1}=\left.\Orange{x_1}{x_2}\right.
\mr{\ and\ }
t_{2}=\left.\Apple{x_3}{x_4}\right. ) \]
and gives an $X$-labelled tree by substituting at the leaves (here,
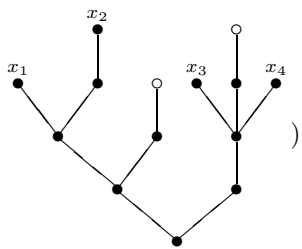
\[\left.\MixedFruit{x_1}{x_2}{x_3}{x_4}\right. ).\]

\item
On the category \Cat\ of small categories and functors, there is the
free strict monoidal category monad. Both \Cat\ and the monad are
cartesian.

\item
A \emph{globular set} is a diagram
\begin{diagram}
\cdots &\pile{\rTo \\ \rTo} &X_n &\pile{\rTo^s \\ \rTo_t} &X_{n-1}
&\pile{\rTo^s \\ \rTo_t} &\cdots &\pile{\rTo \\ \rTo} &X_1
&\pile{\rTo^s \\ \rTo_t} &X_0 \\
\end{diagram}
in \Sets\ satisfying the `globularity equations' $ss=st$ and $ts=tt$.  The
underlying graph of a strict $\omega$-category is a globular set: $X_n$ is
the set of $n$-cells, and $s$ and $t$ are the source and target
functions. One can construct the free strict $\omega$-category monad on the
category of globular sets and show that it is cartesian. Moreover, the
category of globular sets is cartesian, being a presheaf category.

\item
A \emph{double category} may be defined as a category object in
\fcat{Cat}. More descriptively, the graph structure consists of collections
of
\begin{itemize}
	\item 0-cells $A$
	\item horizontal 1-cells $f$
	\item vertical 1-cells $p$
	\item 2-cells $\alpha$
\end{itemize}
and various source and target functions, as illustrated by the picture
\[
\begin{diagram}
A_1		&\rTo^{f_1}		&A_2		\\
\dTo<{p_1}	&\Downarrow \alpha	&\dTo>{p_2}	\\
A_3		&\rTo_{f_2}		&A_4		\\
\end{diagram}
\ .
\]
The category structure consists of identities and composition functions for
2-cells and both kinds of 1-cell, obeying strict associativity, identity and
interchange laws; see \cite{KS} for more details.

More generally, let us define \emph{$n$-cubical set} for any $n\elt\nat$; the
intention is that a 2-cubical set will be the underlying graph of a double
category. So, let \fcat{Cube_n} be the category with
\begin{description}
	\item[objects:] subsets $D$ of $\{ \range{0,1}{n-1} \}$ 
	\item[map $D \go D'$:] the inclusion $D\sub D'$, together with 
	a function $D'\without D \go \{0,1\}$ 
	\item[composition:] place functions side-by-side.
\end{description}
Then we define an $n$-cubical set to be a presheaf on \fcat{Cube_n}. For
instance, we may think of a 2-cubical set $X$ as:
\begin{itemize}
	\item $X\emptyset = \{ \mbox{0-cells} \}$
	\item $X\{0\} = \{ \mbox{horizontal 1-cells} \}$
	\item $X\{1\} = \{ \mbox{vertical 1-cells} \}$
	\item $X\{0,1\} = \{ \mbox{2-cells} \}$
\end{itemize}
and, for instance, the map $\{1\} \go \{0,1\}$ given by
\[
\{0,1\} \without \{1\} = \{0\} \goby{0} \{0,1\}
\]
sends $\alpha\elt X \{0,1\}$ to $p_{1} \elt X\{1\}$, in the diagram above. In
the context of functions $D'\without D \go \{0,1\}$, 0 should be read as
`source' and 1 as `target'.

We may now define a \emph{(strict) $n$-tuple category} to be an $n$-cubical
set together with various compositions and identities, as for double
categories, all obeying strict laws. The category of $n$-cubical sets has on
it the free strict $n$-tuple category monad; both category and monad are
cartesian.

\end{eg}


\clearpage


\section{Multicategories}	\label{sec:multicats}

We now describe what an \Cartpr-multicategory is, where \ust\ is a cartesian
monad on a cartesian category \ess. As mentioned in the Introduction, this
description is a generalization of the (well-known) description of a small
category as a monad object in the bicategory of spans.

We will use the phrase `\Cartpr\ is cartesian' to mean that \ess\ is a
cartesian category and \Mnd\ is a cartesian monad on \ess.

\begin{constn} \label{constn:bicat} \end{constn}
Let \Cartpr\ be cartesian. We construct a bicategory \cat{B} from
\Cartpr, which in the case $\ust=\id$ is the bicategory of spans in
\ess.
\begin{description}
\item[0-cell:] Object $S$ of \ess.
\item[1-cell $R \go S$:] Diagram
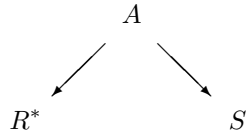
\begin{slopeydiag}
	&	&A	&	&	\\
	&\ldTo	&	&\rdTo	&	\\
R^{*}	&	&	&	&S	\\
\end{slopeydiag}
in \ess.
\item[2-cell $A \go A'$:] Commutative diagram
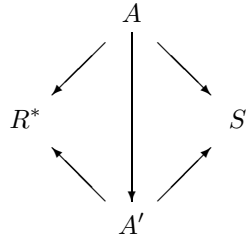
\begin{slopeydiag}
	&	&A	&	&	\\
	&\ldTo	&	&\rdTo	&	\\
R^{*}	&	&\dTo	&	&S	\\
	&\luTo	&	&\ruTo	&	\\
	&	&A'	&	&	\\
\end{slopeydiag}
in \ess.
\item[1-cell composition:] To define this we need to choose particular
pullbacks in \ess, and in everything that follows we assume this has
been done. Take 
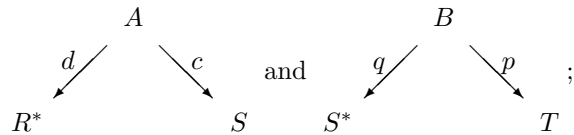
\[\left.
\begin{slopeydiag}
	&	&A	&	&	\\
	&\ldTo<{d}&	&\rdTo>{c}&	\\
R^{*}	&	&	&	&S	\\
\end{slopeydiag}
\right.
\mr{\ and\ }
\left.
\begin{slopeydiag}
	&	&B	&	&	\\
	&\ldTo<{q}&	&\rdTo>{p}&	\\
S^{*}	&	&	&	&T	\\
\end{slopeydiag}
\right.
;
\]
then their composite is given by the diagram
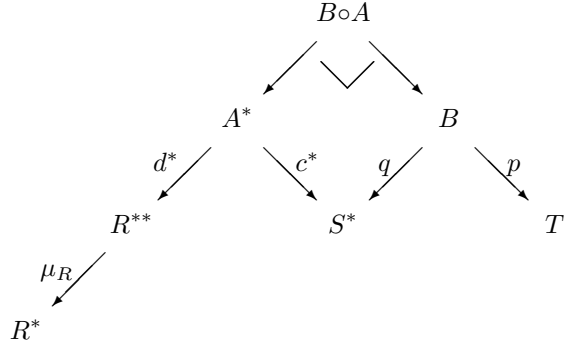
\begin{slopeydiag}
   &       &   &       &   &       &B\of A\Spbk&  &   &       &   \\
   &       &   &       &   &\ldTo  &      &\rdTo  &   &       &   \\
   &       &   &       &A^*&       &      &       &B  &       &   \\
   &       &   &\ldTo<{d^*}&&\rdTo>{c^*}& &\ldTo<{q}& &\rdTo>{p}& \\
   &       &R^{**}&    &   &       &S^*   &       &   &       &T  \\
   &\ldTo<{\mu_R}&&    &   &       &      &       &   &       &   \\
R^*&       &   &       &   &       &      &       &   &       &   \\
\end{slopeydiag}
where the right-angle mark in the top square indicates that the square
is a pullback.

\item[1-cell identities:] The identity on $S$ is
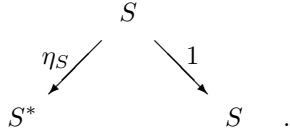
\begin{slopeydiag}
	&	&S	&	&	&\\
	&\ldTo<{\eta_S}&&\rdTo>{1}&	&\\
S^{*}	&	&	&	&S	&.\\
\end{slopeydiag}

\item[2-cell identities and compositions:] Identities and vertical
composition are as in \ess. Horizontal composition is given in an
obvious way.

\end{description}

Because the choice of pullbacks is arbitrary, 1-cell composition does
not obey strict associative and identity laws. That it obeys them up
to invertible 2-cells is a consequence of the fact that \Mnd\ is
cartesian. \done

\begin{defn}
A \emph{monad in a bicategory} consists of a 0-cell $S$, a 1-cell $S
\goby{t} S$, and 2-cells
\begin{diagram}
S&
\pile{\rTo^1	\\ \Downarrow\eta	\\ \rTo~t	
		\\ \Uparrow\mu		\\ \rTo_{t \of t}}
&S,\\
\end{diagram}
such that the diagrams
\[\left.
\begin{diagram}
t\of 1	&\rTo^{t\eta}	&t\of t		&\lTo^{\eta t}	&1\of t	\\
   	&\rdTo<{\diso}	&\dTo<{\mu}	&\ldTo>{\diso}	&	\\
   	&     		&t  		&     		&   	\\
\end{diagram}
\right.
\ \ \ 
\left.
\begin{diagram}
   		&     		&t\of(t\of t)	&\rTo^{t\mu}	&t\of t\\
   		&\ldLine<{\diso}&   		&     		&   \\
(t\of t)\of t	&     		&   		& 		&\dTo>{\mu}\\
\dTo<{\mu t}	&     		&   		&     		&   \\
t\of t		&     		&\rTo_{\mu}	&     		&t  \\
\end{diagram}
\right.\]
commute.
\end{defn}

\begin{defn}
Let \Cartpr\ be cartesian. Then an \emph{\Cartpr-multicategory} is a
monad in the associated bicategory \cat{B} of
Construction~\ref{constn:bicat}.
\end{defn}

An \Cartpr-multicategory therefore consists of a diagram
\spaan{A}{S^*}{S}{d}{c} in \ess\ and maps $S \goby{\ids} A$, $A\of A
\goby{\comp} A$ satisfying associative and identity laws. Think of $S$
as `objects', $A$ as `arrows', $d$ as `domain' and $c$ as `codomain'.
Such an $A$ will be called an \Cartpr-multicategory \emph{on $S$}, or
if $S=1$ an \emph{\Cartpr-operad}. (Baez and Dolan, in \cite{BD}, use
`operad' or `typed operad' for the same kind of purpose as we use
`multicategory', and `untyped operad' where we use `operad'.)

It is inherent that everything is small: when $\ess=\Sets$, for instance, the
objects and arrows form sets, not classes. Dealing with large multicategories
instead does not appear to present any problem in practice.  If we wanted to
perform category theory enriched over a multicategory (in the place of a
monoidal category), then the use of large multicategories would be
necessary. (See~\ref{fds:struc} for further remarks on the relationship
between multicategories and monoidal categories.)

\begin{defn}
Let \Cartpr\ be cartesian.
\begin{enumerate}
\item \label{defn:graph}
An \emph{\Cartpr-graph (on $S$)} is a diagram \gph{A}{S} in
\ess. A \emph{map of \Cartpr-graphs}
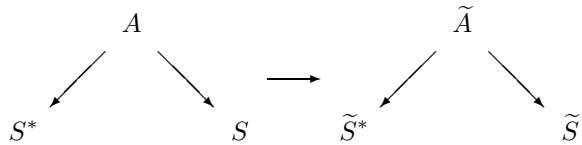
\[\left.
\begin{slopeydiag}
	&	&A	&	&	\\
	&\ldTo	&	&\rdTo	&	\\	
S^{*}	&	&	&	&S	\\
\end{slopeydiag}
\right.
\go
\left.
\begin{slopeydiag}
		&	&\twid{A}&	&	\\	
		&\ldTo	&	&\rdTo	&	\\
\twid{S}^{*}	&	&	&	&\twid{S}\\
\end{slopeydiag}
\right.\]
is a pair \pr{A \goby{f} \twid{A}}{S \goby{g} \twid{S}} of maps in
\ess\ such that
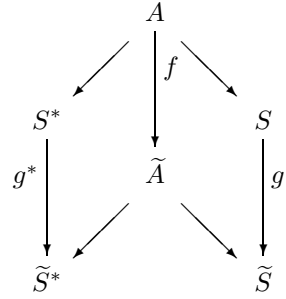
\begin{slopeydiag}
	&	&A	&	&	\\
	&\ldTo	&\dTo>{f}&\rdTo	&	\\
S^{*}	&	&	&	&S	\\
\dTo<{g^*}&	&\twid{A}&	&\dTo>{g}\\
	&\ldTo	&	&\rdTo	&	\\
\twid{S}^{*}&	&	&	&\twid{S}\\
\end{slopeydiag}
commutes.

\item A \emph{map of \Cartpr-multicategories} $A \go \twid{A}$ (with graphs
as in~\ref{defn:graph}) is a map \pr{f}{g} of their graphs such that the
diagrams 
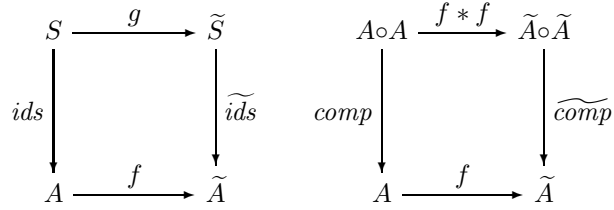
\[\left.
\begin{diagram}
S		&\rTo^{g}		&\twid{S}		\\	
\dTo<{\ids}	&			&\dTo>{\twid{\ids}}	\\
A		&\rTo^{f}		&\twid{A}		\\
\end{diagram}
\right.
\ \ \ 
\left.
\begin{diagram}
&A\of A		&\rTo^{f*f}		&\twid{A}\of\twid{A}	\\
&\dTo<{\comp}	&			&\dTo>{\twid{\comp}}	\\
&A		&\rTo^{f}		&\twid{A}		\\	
\end{diagram}
\right.\]
commute.

\end{enumerate}
\end{defn}

\begin{remarks}{rmks:maps}
\item The map $A\of A \goby{f*f} \twid{A}\of\twid{A}$  just mentioned
is the horizontal composite of 2-cells in the bicategory \cat{B} of
Construction~\ref{constn:bicat}. That is, $f*f$ is the unique map in
\ess\ making
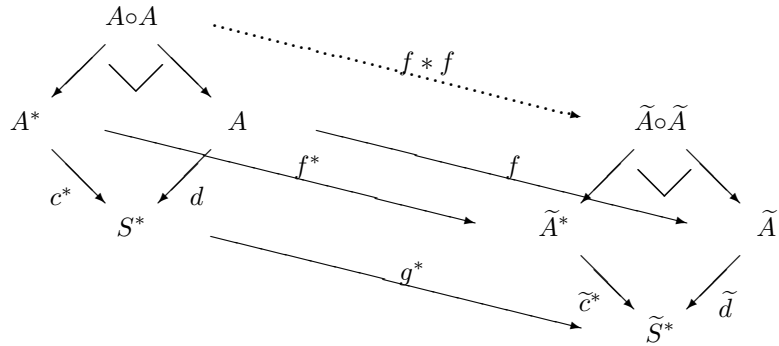
\begin{slopeydiag}
   &   &A\of A\Spbk&   &   &
   &   &   &   &   & 
   &   &   &   &   \\
   &\ldTo&   &\rdTo\rdGet(10,2)^{f*f}&   &
   &   &   &   &   &
   &   &   &   &   \\
A^*&   &   &   &A  &
   &   &   &   &   &
   &   &\twid{A}\of\twid{A}\Spbk&   &   \\
   &\rdTo<{c^*}\rdTo(10,2)^{f^*}&   &\ldTo>{d}&   &
\rdTo(10,2)^{f}&   &   &   &   &
   &\ldTo&   &\rdTo&   \\
   &   &S^*&   &   &
   &   &   &   &   &
\twid{A}^*&   &   &   &\twid{A}\\
   &   &   &\rdTo(10,2)^{g^*}&     &
   &   &   &   &   &
   &\rdTo<{\twid{c}^*}&   &\ldTo>{\twid{d}}&   \\
   &   &   &   &   &
   &   &   &   &   &
   &   &\twid{S}^*&   &   \\
\end{slopeydiag}
commute.

\item Fix $S\elt\ess$. Then we may consider the category of
\Cartpr-graphs on $S$, whose morphisms \pr{A \goby{f} \twid{A}}
{S \goby{g} S} all have $g=1$. This is just the slice category
$\frac{\ess}{S^{*}\times S}$. It is also the full sub-bicategory of
\cat{B} whose only object is $S$, and is therefore a monoidal
category. The category of \Cartpr-multicategories on $S$ is then the
category $\Mon(\frac{\ess}{S^{*}\times S})$ of monoids in 
$\frac{\ess}{S^{*}\times S}$.

\item A choice of pullbacks in \ess\ was made; changing that choice
gives an isomorphic category of \Cartpr-multicategories. 

\end{remarks}
 
\begin{eg}{eg:multicats}

\item
Let $\Cartpr=\pr{\Sets}{\id}$. Then \cat{B} is the
bicategory of spans, and a monad in \cat{B} is just a (small)
category. Thus categories are \pr{\Sets}{\id}-multicategories.
Functors are maps of such.

\item
Let $\Cartpr=\pr{\Sets}{\mr{free\ monoid}}$. Specifying an \Cartpr-graph
\graph{A}{S}{d}{c} is equivalent to specifying a set
\mtihom{A}{s_1}{s_n}{s} for each \range{s_1}{s_n,s}\elt$S$ ($n\geq
0$); if $a\elt\mtihom{A}{s_1}{s_n}{s}$ then we write
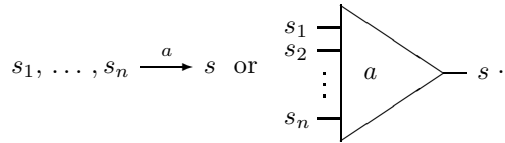
\[
\range{s_1}{s_n} \goby{a} s
\mr{\ \ or\ \ }
\ctransistor{a}{s_1}{s_2}{s_n}{s}
\mr{\ .}
\]

In the associated bicategory, the identity 1-cell
\graph{S}{S}{\eta_S}{1} on $S$ has
\[
\mtihom{S}{s_1}{s_n}{s}=
\left\{
\begin{array}{ll}
1		&\mr{if\ }n=1\mr{\ and\ }s_{1}=s	\\
\emptyset	&\mr{otherwise.}
\end{array}
\right.
\]
The composite 1-cell $A\of A$ is
\[
\{\pr{\bftuple{a_1}{a_n}}{a} \such da=\bftuple{ca_1}{ca_n}\},
\]
i.e.\ is the set of diagrams
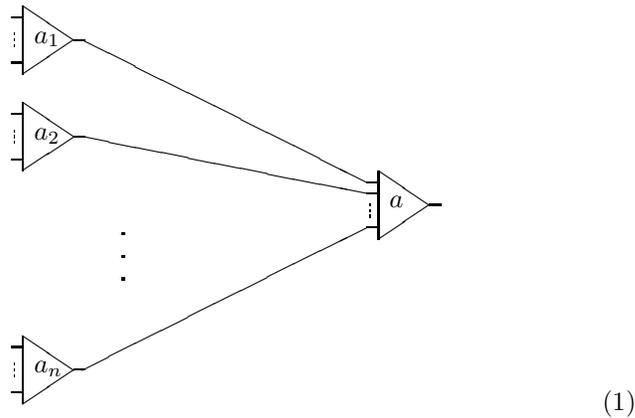
\begin{equation} \label{diag:comptrans}
\comptrans{a_1}{a_2}{a_n}{a}
\end{equation}
with the evident domain and codomain functions.

We then have a function \ids\ assigning to each $s\elt S$ a member of
\mtiendo{A}{s}, and a function \comp\ composing diagrams
like~(\ref{diag:comptrans}). These are required to obey associative and
identity laws. Thus a \pr{\Sets}{\mr{free\ monoid}}-multicategory is
just an `ordinary' non-symmetric multicategory. A \pr{\Sets}{\mr{free\
monoid}}-operad is a non-symmetric May operad (on \Sets).

\item
One should not conclude from
Example~\ref{eg:cart-mnds}(\ref{eg:comm-not-cart}) that it is
impossible in our system to describe the symmetric operads of \cite{May} or
\cite{BD}. The reason why not lies in the difference between having an
isomorphism 
\[
\mtihom{A}{s_1}{s_n}{s} \iso \mtihom{A}{s_{\sigma 1}}{s_{\sigma n}}{s}
\]
for each permutation $\sigma$, and having actual equalities.

\item
Let $\ess=\Sets$, and consider the exceptions monad $\dashbk+1$
of~\ref{eg:cart-mnds}(\ref{eg:exceptions-mnd}). A
\pr{\Sets}{$\dashbk+1$}-graph is a diagram \spaan{A}{S+1}{S}{d}{c} of
sets; this is like an ordinary \pr{\Sets}{\id}-graph on $S$, except
that some arrows have domain 0---an extra element not in $S$. (Thus
$1=\{0\}$ here.) If we set
\[
Y(s) = \{a\elt A \such da=0\}
\]
for each $s\elt S$, then a multicategory structure on the graph provides a
function
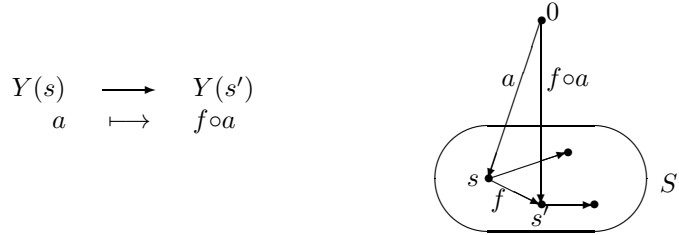
\[
\begin{array}{rcl}
Y(s)	&\go	&Y(s')	\\
a	&\goesto&f\of a
\end{array}
\cdiscfibpic
\]
for each $f\elt A$ with $d(f)=s\elt S$ and $c(f)=s'$. It also provides
a category structure on \spaan{C}{S}{S}{d}{c}, where $C=\{a\elt A\such
da\elt S\}$. Thus a \pr{\Sets}{\dashbk+1}-multicategory turns out to
be just a (small) category \scat{C} together with a functor
$\scat{C}\go\Sets$. (Similarly, a
\pr{\Sets}{\dashbk+E}-multicategory is a category \scat{C} together
with an $E$-indexed family of functors $\scat{C}\go\Sets$.)

Another way to put this is that a \pr{\Sets}{\dashbk+1}-multicategory is a
discrete fibration (between small categories, where the codomain here is
$\scat{C}^{\op}$). In fact, the category of
\pr{\Sets}{\dashbk+1}-multicategories is the category of discrete fibrations.

\item
Let $M$ be a monoid and $\Cartpr = \pr{\Sets}{M\times\dashbk}$. Then an
\Cartpr-multicategory consists of a category \scat{C} together with a functor
$\scat{C}\go M$.

\item
Let $\Cartpr = \pr{\Sets}{\mr{tree\ monad}}$, as
in~\ref{eg:cart-mnds}(\ref{eg:tree-mnd}). An \Cartpr-multicategory
consists of a set $S$ of objects, and sets like
\[
A\left(\Pear{s_1}{s_2},s\right)
\]
($s_{1}, s_{2}, s \elt S$), together with a unit element of each
$A(\nl{s},s)$ and composition functions like
\begin{eqnarray*}
\left\{	A\left( \Orange{r_1}{r_2},s_1\right)
	\times
	A\left( \Apple{r_3}{r_4},s_2\right)
\right\}
\times
A\left( \Pear{s_1}{s_2},s\right)\\
\go
A\left( \MixedFruit{r_1}{r_2}{r_3}{r_4},s\right)
\end{eqnarray*}
($r_{1}, r_{2}, r_{3}, r_{4} \elt S$). These are to satisfy associativity
and identity laws.

When $S=1$, so that we're considering \Cartpr-operads, the graph structure is
comprised of sets like $A\left(\Pear{}{}\right)$.

The \Cartpr-multicategories are a simpler version of Soibelman's
pseudo-monoidal categories (\cite{Soi}); they omit the aspect of maps between
trees. A similar relation is borne to Borcherds' relaxed multilinear
categories (\cite{Bor}).

\item
When $\ess=\fcat{Cat}$ and \ust\ is the free strict monoidal category monad,
an \Cartpr-operad is what Soibelman calls a strict monoidal 2-operad in
\cite[2.1]{Soi}.

\item
Let $\Cartpr=\pr{\mr{Globular\ sets}}{\mr{free\ strict\
}\omega\mr{-category}}$. An \Cartpr-operad is exactly what Batanin
calls an operad (or `an $\omega$-operad in \emph{Span}'; see~\cite{Bat}).

\item
Operads for $\Cartpr = \pr{n\mr{-cubical\ sets}}{\mr{free\ strict\
}n\mr{-tuple\ category}}$ can be understood in much the same way as Batanin's
operads. For instance, a cell in the free strict $\omega$-category on the
terminal globular set can be represented as a tree\footnote{These are not the
same kind of trees as in \ref{eg:cart-mnds}(\ref{eg:tree-mnd}); see
\cite{Bat}}; a cell in the free strict $n$-tuple category on the terminal
$n$-cubical set can be represented as a cuboid (or the sequence of natural
numbers which are its edges' lengths). A Batanin operad associates to each
tree a set, and has composition functions corresponding to the combining of
trees; a cubical operad associates to each cuboid a set, and has composition
functions corresponding to the combining of cuboids.

\end{eg}


\clearpage


\section{Algebras}	\label{sec:algs}

We want to define a category of algebras for any
\Cartpr-multicategory. In the case $\Cartpr=\pr{\Sets}{\id}$, where
we are dealing with a plain category \scat{C}, the category of
algebras should be \ftrcat{\scat{C}}{\Sets}. We will take inspiration
from the following:

\begin{lemma} \label{lemma:csets-adjn}
Let \scat{C} be a small category, and let $\scat{C}_{0}$ denote the
set of objects of \scat{C} or the discrete category thereon. Then the
forgetful functor 
$\ftrcat{\scat{C}}{\Sets} \go \ftrcat{\scat{C}_0}{\Sets}$
is monadic.
\end{lemma}

\pf\ This is easily verified without use of the adjoint functor or
monadicity theorems. Here, we just describe what the induced monad $T$
does to an object $X$ of $\ftrcat{\scat{C}_0}{\Sets}$: if $C \elt
\scat{C}$ then
\[
(TX)C=\coprod_{\stackrel{D \goby{f} C}{\mr{\ in\ } \scat{C}}} XD.
\]
\done

\begin{discussion} \label{dis:ind-mnd} \end{discussion}
Lemma~\ref{lemma:csets-adjn} describes $\ftrcat{\scat{C}}{\Sets}$ as
the category of algebras for a certain monad $T$ on
$\ftrcat{\scat{C}_0}{\Sets}$. But there is an equivalence
\[
\ftrcat{\scat{C}_0}{\Sets} \eqv \Sets/\scat{C}_{0}
\]
of categories, so we also obtain a monad $T'$ on $\Sets/\scat{C}_{0}$. How
can $T'$ be described directly? Take an object \slob{X}{p}{\scat{C}_0}
of $\Sets/\scat{C}_{0}$: this corresponds to the object $X$ of 
$\ftrcat{\scat{C}_0}{\Sets}$ with $XC=p^{-1}\{C\}$. Now
\begin{eqnarray*}
(TX)C	&=&	\coprod_{D\goby{f}C} XD		\\
	&=&	\{\pr{f}{v}\such f\elt\scat{C}_{1},\ v\elt X(\dom f), 
				 \cod f=C\}\\
	&=&	\{\pr{f}{v}\elt\scat{C}_{1}\times X \such p(v)=\dom f
							,\ \cod f=C\}
\end{eqnarray*}
(where $\scat{C}_{1}=\{$arrows of $\scat{C}\}$), so if
$T'X=\slob{X'}{p'}{\scat{C}_0}$ then $X'$ is the pullback
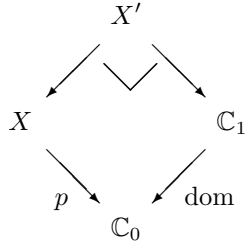
\begin{slopeydiag}
	&	&X'\Spbk&	&	\\
	&\ldTo	&	&\rdTo	&	\\
X	&	&	&	&\scat{C}_{1}	\\
	&\rdTo<{p}&	&\ldTo>{\dom}&	\\
	&	&\scat{C}_{0}&	&	\\
\end{slopeydiag}
and $p'(f,v)=\cod f$. Thus $T'X$ is the right-hand diagonal of the diagram
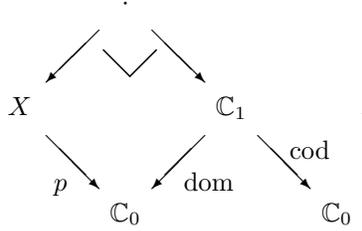
\[\left.
\begin{slopeydiag}
	&	&\cdot\Spbk&	&	&	&	\\
	&\ldTo	&	&\rdTo	&	&	&	\\
X	&	&	&	&\scat{C}_{1}&	&	\\
	&\rdTo<{p}&	&\ldTo>{\dom}&	&\rdTo>{\cod}&	\\
	&	&\scat{C}_{0}&	&	&	&\scat{C}_{0}\\
\end{slopeydiag}
\right.
.
\]
We are now ready to generalize to any cartesian \Cartpr.

\begin{constn} \label{constn:ind-mnd} \end{constn}
Let \Cartpr\ be cartesian and $S\elt\ess$: then any \Cartpr-multicategory on
$S$ gives rise to a monad on $\ess/S$.

Let \graph{A}{S}{d}{c} be the multicategory. We describe a monad \blbk\ on
\Slice; in what follows, we'll write $\bktdslob{X}{p}{S}^{\blob} =
\slob{X_\blob}{p_\blob}{S}$, etc.
\begin{itemize}
\item \slob{X_\blob}{p_\blob}{S} is the right-hand diagonal of the diagram
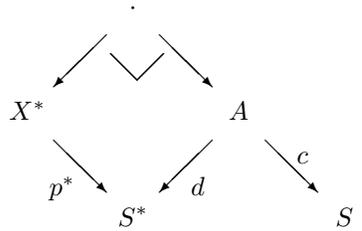
\[\left.
\begin{slopeydiag}
	&	&\cdot\Spbk&	&	&	&	\\
	&\ldTo	&	&\rdTo	&	&	&	\\
X^*	&	&	&	&A	&	&	\\
	&\rdTo<{p^*}&	&\ldTo>{d}&	&\rdTo>{c}&	\\
	&	&S^*	&	&	&	&S	\\
\end{slopeydiag}
\right.
.
\]

\item If
$\left.
\begin{triangdiag}
X	&	&\rTo^{f}&	&Y	\\
	&\rdTo<{p}&	&\ldTo>{q}&	\\
	&	&S	&	&	\\
\end{triangdiag}
\right.$
is a map in \Slice, then 
$\left.
\begin{triangdiag}
X_\blob	&		&\rTo^{f^\blob}	&		&Y_\blob\\
	&\rdTo<{p_\blob}&		&\ldTo>{q_\blob}&	\\
	&		&S		&		&	\\
\end{triangdiag}
\right.$
is the unique map making
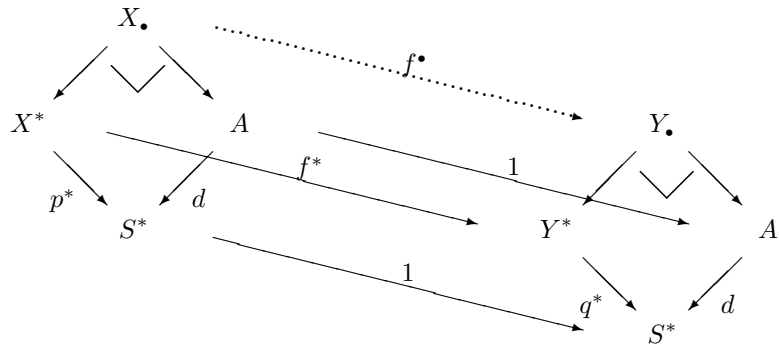
\begin{slopeydiag}
   &   &X_{\blob}\Spbk&   &   &
   &   &   &   &   & 
   &   &   &   &   \\
   &\ldTo&   &\rdTo\rdGet(10,2)^{f^\blob}&   &
   &   &   &   &   &
   &   &   &   &   \\
X^*&   &   &   &A  &
   &   &   &   &   &
   &   &Y_{\blob}\Spbk&   &   \\
   &\rdTo<{p^*}\rdTo(10,2)^{f^*}&   &\ldTo>{d}&   &
\rdTo(10,2)^{1}&   &   &   &   &
   &\ldTo&   &\rdTo&   \\
   &   &S^*&   &   &
   &   &   &   &   &
Y^*&   &   &   &A\\
   &   &   &\rdTo(10,2)^{1}&     &
   &   &   &   &   &
   &\rdTo<{q^*}&   &\ldTo>{d}&   \\
   &   &   &   &   &
   &   &   &   &   &
   &   &S^*&   &   \\
\end{slopeydiag}
commute.

\item The unit at \slob{X}{p}{S} is given by
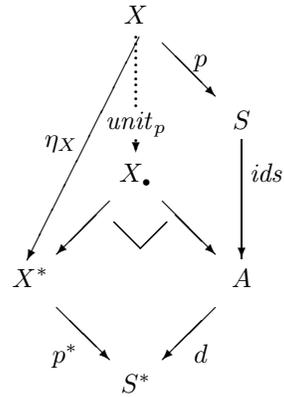
\[\left.
\begin{slopeydiag}
	&			&X		&	&	\\
	&\ldTo(2,5)<{\eta_X}	&		&\rdTo>{p}&	\\
	&			&\dGet~{\unit_p}	&	&S	\\
	&			&X_{\blob}\Spbk	&	&\dTo>{\ids}\\
	&\ldTo			&		&\rdTo	&	\\
X^{*}	&			&		&	&A	\\
	&\rdTo<{p^*}		&		&\ldTo>{d}&	\\
	&			&S^{*}		&	&	\\
\end{slopeydiag}
\right.
\ \ \ \ \ \ .
\]

\item For multiplication, we have a commutative diagram
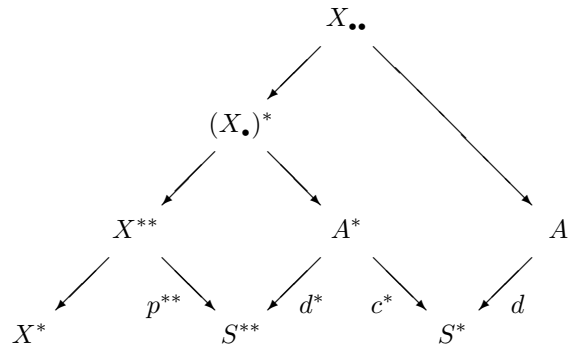
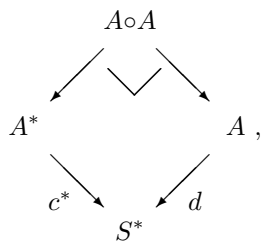
\begin{slopeydiag}
   &	&	&	&	&	&X_{\blob\blob}&&    &	&  \\
   &	&	&	&	&\ldTo	&	&\rdTo(4,4)& &	&  \\
   &	&	&	&(X_{\blob})^*&	&	&	&    &	&  \\
   &	&	&\ldTo	&	&\rdTo	&	&	&    &	&  \\
   &	&X^{**}	&	&	&	&A^*	&	&    &	&A \\
   &\ldTo&	&\rdTo<{p^{**}}&&\ldTo>{d^*}&	&\rdTo<{c^*}&&\ldTo>{d}&\\
X^*&	&	&	&S^{**}	&	&	&	&S^* &  &  \\
\end{slopeydiag}
and a pullback square
\[
\begin{slopeydiag}
	&	&A\of A\Spbk&	&	\\
	&\ldTo	&	&\rdTo	&	\\
A^{*}	&	&	&	&A	\\
	&\rdTo<{c^*}&	&\ldTo>{d}&	\\
	&	&S^*	&	&	\\
\end{slopeydiag}
\mr{\ ,}
\]
and so in particular there are maps
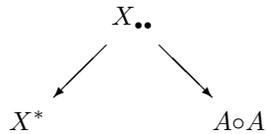
\[\left.
\begin{slopeydiag}
	&	&X_{\blob\blob}	&	&	\\
	&\ldTo	&		&\rdTo	&	\\
X^{*}	&	&		&	&A\of A	\\
\end{slopeydiag}
\right.
.
\]
The multiplication at \slob{X}{p}{S} is given by
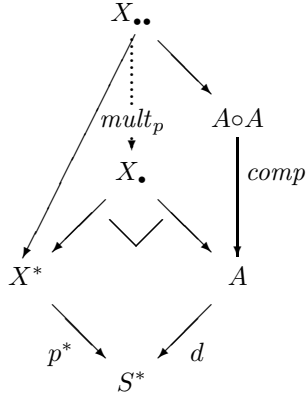
\[	
\begin{slopeydiag}
	&			&X_{\blob\blob}	&	&	\\
	&\ldTo(2,5)		&		&\rdTo	&	\\
	&			&\dGet~{\mult_p}&	&A\of A	\\
	&			&X_{\blob}\Spbk	&	&\dTo>{\comp}\\
	&\ldTo			&		&\rdTo	&	\\
X^{*}	&			&		&	&A	\\
	&\rdTo<{p^*}		&		&\ldTo>{d}&	\\
	&			&S^{*}		&	&	\\
\end{slopeydiag}
\ \ \ \ \ \ \ \ \ \ \ .
\]
\end{itemize}

It is now straightforward to check that \Imnd\ forms a monad on \Slice, and
that when $\Cartpr=\pr{\Sets}{\id}$ this is the monad of
Discussion~\ref{dis:ind-mnd}. \done

\begin{defn}
Let \Cartpr\ be cartesian and \gph{A}{S} an \Cartpr-multicategory on
$S\elt\ess$. Then the \emph{category of algebras} for the multicategory,
$\Alg(A)$, is the category of algebras for the associated monad on \Slice.
\end{defn}

With the \pr{\Sets}{\id} case of plain categories in mind, we would expect a
map $A \go A'$ of multicategories to yield a functor $\Alg(A) \og \Alg(A')$.
This is indeed the case; moreover, one can define a \emph{transformation}
between two maps of multicategories, and such a transformation leads to a
natural transformation
\[
\begin{diagram}
\Alg(A)	&\pile{\lTo	\\ \Downarrow	\\ \lTo}	&\Alg(A')\\
\end{diagram}
\ .
\]
The question also arises of `Kan extensions': left and right adjoints to the
functor $\Alg(A) \og \Alg(A')$. These issues will not be discussed any
further here.

\begin{eg}{eg:algs}

\item 
When $\Cartpr=\pr{\Sets}{\id}$, $\Alg(\scat{C}) \iso
\ftrcat{\scat{C}}{\Sets}$.

\item
When $\Cartpr=\pr{\Sets}{\mr{free\ monoid}}$, so that an
\Cartpr-multicategory is a multicategory of the familiar kind, we already
have an idea of what an algebra for $A$ should be: a `multifunctor
$A\go\Sets$'. That is, an algebra for $A$ should consist of:
	\begin{itemize}
	\item for each $s\elt S$, a set $X(s)$
	\item for each $\range{s_1}{s_n}\goby{f}s$ in $A$, a function
	$X(s_{1})\times\cdots\times X(s_{n}) \go X(s)$, preserving identities
	and composition.  
	\end{itemize}
In fact, this is the same as the definition of algebra just given. One way to
see this is to work through Lemma~\ref{lemma:csets-adjn} and
Discussion~\ref{dis:ind-mnd}, changing `category' to `(ordinary)
multicategory': \ftrcat{A}{\Sets} is monadic over \ftrcat{S}{\Sets}, 
the monad $T$ being given by
\[
(TX)s = \coprod_{\range{s_1}{s_n}\goby{f}s} X(s_{1})\times\cdots\times
X(s_{n}).
\]
Alternatively, one can calculate directly: if \slob{X}{p}{S} and we put
$X(s)=p^{-1}\{s\}$ then
\begin{eqnarray*}
X_{\blob}	&=&	\{\pr{\bftuple{x_1}{x_n}}{f} \such
			df=\bftuple{px_1}{px_n}\}\\
		&=&	\{X(s_{1})\times\cdots\times X(s_{n}) \times
			\mtihom{A}{s_1}{s_n}{s} \such
			\range{s_1}{s_{n},s}\elt S\},
\end{eqnarray*}
and an algebra structure on \slob{X}{p}{S} therefore consists of a function
\[
X(s_{1})\times\cdots\times X(s_{n}) \go X(s)
\]
for each member of \mtihom{A}{s_1}{s_n}{s}, subject to certain laws.

\item
When $\Cartpr=\pr{\Sets}{\dashbk+1}$, an \Cartpr-multicategory is an ordinary
category \scat{C} together with a functor $\scat{C}\goby{Y}\Sets$. A
\pr{\scat{C}}{Y}-algebra is then a functor $\scat{C}\go\Sets$ together with a
natural transformation
\[
\begin{diagram}
\scat{C}	&\pile{\rTo^Y \\ \Downarrow \\ \rTo}	&\Sets	\\
\end{diagram}
\mr{\ .}
\] 
In terms of fibrations, an \Cartpr-multicategory is a discrete fibration
$Y$ over a small category \scat{B} (=$\scat{C}^{\op}$), and an algebra for
$Y$ consists of another discrete fibration $X$ over \scat{B} together with a
map from $Y$ to $X$ (of fibrations over \scat{B}).

\item
Let \Cartpr\ be the tree monad on \Sets; for simplicity, let us just consider
\emph{operads} $A$ for \Cartpr---thus the object-set $S$ is 1. An algebra for
$A$ consists of a set $X$ together with a function $X_{\blob}\go X$
satisfying some axioms. One can calculate that an element of $X_{\blob}$
consists of an $X$-labelling of a tree $T$ together with a member of $A(T)$.
An $X$-labelling of an $n$-leafed tree $T$ is just a member of $X^n$, so one
can view the algebra structure $X_{\blob}\go X$ on $X$ as: for each number
$n$, $n$-leafed tree $T$, and element of $A(T)$, a function $X^{n}\go{X}$.
These functions are required to be compatible with glueing of trees in an
evident way.

\item
For $\Cartpr=\pr{\mr{Globular\ sets}}{\mr{free\ strict\ }
\omega\mr{-category}}$, Batanin constructs a certain operad $K$, the
`universal contractible operad' (see \cite{Bat}). He then defines a weak
$\omega$-category to be an algebra for $K$.

\item	\label{eg:terminal}
The graph \graph{1^*}{1}{1}{!} is terminal amongst all \Cartpr-graphs. It
carries a unique multicategory structure, since a terminal object in a
monoidal category always carries a unique monoid structure. It then becomes
the terminal \Cartpr-multicategory. The induced monad on $\ess/1$ is just
\Mnd, and so an algebra for the terminal multicategory is just an algebra for
\ust. This can aid recognition of when a theory of operads or
multicategories fits into our scheme. For instance, if we were to read
Batanin's paper and learn that, in his terminology, an algebra for the
terminal operad is a strict $\omega$-category (\cite[p.\ 51, example
3]{Bat}), then we might suspect that his operads were \Cartpr-operads for the
free strict $\omega$-category monad \ust\ on some suitable category \ess---as
indeed they are.

\end{eg}


\clearpage


\section{Further Developments}	\label{sec:fds}

We finish with a collection of loose topics. Some of them have not been
worked out in full detail; others have, but their relevance is unclear.

Section~\ref{fds:free} is a brief explanation of the process of forming the
free multicategory on a graph, and allows descriptions of both the set of
opetopes and, for any cartesian \Cartpr, a multicategory whose algebras are
the \Cartpr-multicategories. Section~\ref{fds:slicing} explains `slicing', a
generalization of the Grothendieck construction and another important
component of the Baez-Dolan theory. In~\ref{fds:struc} we discuss the
relationship between multicategories and monoidal
categories;~\ref{fds:multicat} throws more light on this by describing the
ways in which the assignment
\[
\Cartpr	\goesto	\Multicatof{\Cartpr}
\]
is functorial. Next, in~\ref{fds:endo}, we associate to any object of \Slice\
an \Cartpr-multicategory on S; this relates to the usual definition of
`algebra' for the operads of May, Baez-Dolan and Batanin.
Sections~\ref{fds:characterization} and~\ref{fds:bicat} each provide an
alternative description of what an \Cartpr-multicategory is, one in terms of
monads and the other in terms of bicategories.

Each section can be read independently of the others. Throughout, we will
denote by \Graphof{\Cartpr} and \Multicatof{\Cartpr} the categories of
\Cartpr-graphs and \Cartpr-multicategories. Sections~\ref{fds:multicat}
and~\ref{fds:characterization} also need the definitions in the following
paragraphs; the rest do not.

\label{p:st-defns}%
Let $T$ and $\twid{T}$ be monads on respective categories \cat{C} and
\twid{\cat{C}}. Then a \emph{monad functor} $\pr{\cat{C}}{T}
\goby{\pr{P}{\phi}} \pr{\twid{\cat{C}}}{\twid{T}}$ consists of a functor
$\cat{C} \goby{P} \twid{\cat{C}}$ together with a natural transformation
\begin{ntdiag}
\cat{C}		&	&\rTo^{T}	&	&\cat{C}	\\
		&	&		&\	&		\\
\dTo<{P}	&	&\ruTo>{\phi}	&	&\dTo>{P}	\\
		&\	&		&	&		\\
\twid{\cat{C}}	&	&\rTo_{\twid{T}}&	&\twid{\cat{C}}	\\
\end{ntdiag}
such that
\begin{diagram}
\twid{T}^{2}P	&\rTo^{\twid{T}\phi}&\twid{T}PT	&\rTo^{\phi T}&PT^2	\\
\dTo<{\twid{\mu}P}&		&		&	&\dTo>{P\mu}	\\
\twid{T}P	&		&\rTo_{\phi}	&	&PT		\\
\end{diagram}
and a similar diagram involving units commute. If $\pr{\cat{C}}{T}
\goby{\pr{Q}{\psi}} \pr{\twid{\cat{C}}}{\twid{T}}$ is another monad functor
then a \emph{monad functor transformation} $\pr{P}{\phi} \go \pr{Q}{\psi}$ is
a natural transformation $P \goby{\alpha} Q$ such that $(\alpha T)\of\phi =
\psi\of(\twid{T}\alpha)$. There is consequently a 2-category \fcat{Mnd},
whose 0-cells are pairs \pr{\cat{C}}{T}, whose 1-cells are monad functors,
and whose 2-cells are monad functor transformations.

There is the dual notion of a \emph{monad opfunctor}, which is just like a
monad functor except that $\phi$ travels in the opposite direction;
similarly, \emph{monad opfunctor transformations}. This gives another
2-category, \fcat{Mnd'}. All of these definitions are taken directly from
Street's paper \cite{St}.

\label{p:cart-st-defns}%
A monad opfunctor \pr{P}{\phi} will be called \emph{cartesian} if $P$
preserves pullbacks and $\phi$ is a cartesian natural transformation, but a
monad functor \pr{P}{\phi} will be called \emph{cartesian} just if $P$
preserves pullbacks. (This is an unhappy situation; the reason for these
definitions is that they give the conditions necessary for the constructions
of~\ref{fds:multicat} to work.) We then define \fcat{CartMnd}, the
sub-2-category of \fcat{Mnd} consisting of cartesian pairs
\pr{\cat{C}}{T}, cartesian monad functors, and all monad functor
transformations; similarly, the sub-2-category \fcat{CartMnd'} of
\fcat{Mnd'}.

\subsection{Free Multicategories}	\label{fds:free}

\subsubsection*{The Free Multicategory Functor}

Let \Cartpr\ be cartesian. Subject to certain further conditions on \Cartpr,
which I won't mention except to say that they hold for the opetopic
construction below, the following are true:
\begin{itemize}
	\item the forgetful functor
		\[
		\Multicatof{\Cartpr} \go \Graphof{\Cartpr}
		\]
	has a left adjoint, the `free \Cartpr-multicategory functor'
	\item the adjunction is monadic
	\item the monad on \Graphof{\Cartpr} is cartesian
	\item all of the above statements are also true for the
 	forgetful functor
		\[
		\Cartpr\mbox{-Multicats on }S \go
		\Cartpr\mbox{-Graphs on }S,
		\]
	for any $S \elt \ess$.
\end{itemize}
(It follows from this and the general theory of monads that any multicategory
is a quotient of a free multicategory; this corresponds to the presentation
of a multicategory by generators and relations.)

\subsubsection*{The Multicategory Multicategory}

Take the free \Cartpr-multicategory monad \ush\ on the category $\ess' =
\Graphof{\Cartpr}$ (for suitable \Cartpr). Then \Cartpr-multicategories are
algebras for \ush. By Example~\ref{eg:algs}(\ref{eg:terminal}), this means
that the terminal \pr{\ess'}{\ush}-multicategory has as its algebras the
\Cartpr-multicategories.

(Related to this is the Baez-Dolan construction of the `$S$-operad operad'
for any object-set $S$: an operad whose algebras are the operads on $S$. To
make sense of the last sentence in our language, read `multicategory' for
`operad'.)

\subsubsection*{Opetopes}

The free multicategory functor enables us to construct the sets $S_n$ of
$n$-opetopes ($n\elt\nat$), as developed in \cite{BD}. (See also \cite{Baez}
for a softer account.) Start with $S_{0}=1$ and $T_{0}=\id$; that is, $T_0$
is the identity monad on $\Sets = \Sets/S_{0}$. Now suppose we have a set
$S_n$ and a cartesian monad $T_n$ on $\Sets/S_{n}$. The terminal object of
$\Sets/S_{n}$ is \slob{S_n}{1}{S_n}; write
\[
T_{n}\bktdvslob{S_n}{1}{S_n} = \vslob{S_{n+1}}{}{S_n}.
\]
The category of \pr{\Sets/S_{n}}{T_n}-graphs on \slob{S_n}{1}{S_n} is 
\[
\frac{\Sets/S_{n}}{T_n(S_{n}\go S_{n})}
=
\frac{\Sets/S_{n}}{S_{n+1} \go S_{n}}
\iso
\frac{\Sets}{S_{n+1}},
\]
so the monad `free \pr{\Sets/S_{n}}{T_n}-multicategory on 1', $T_{n+1}$, is a
cartesian monad on the category $\Sets/S_{n+1}$.

This defines the sets $S_n$; let us look at $n=$0, 1 and 2. First of all,
$S_{0} = 1$ and $T_{0} = \id$. Then
\[
\vslob{S_1}{}{S_0}
=
T_{0}\bktdvslob{S_0}{}{S_0},
\]
i.e.\ $S_{1}=1$, and $T_{1}$ is the monad `free
\pr{\Sets/S_{0}}{T_0}-multicategory on 1', i.e.\ the free monoid monad. Next,
\[
\vslob{S_2}{}{S_1}
=
T_{1}\bktdvslob{S_1}{}{S_1},
\]
i.e.\ $S_2$ is the free monoid \nat\ on the set $S_{1}=1$; $T_2$ is the monad
`free \pr{\Sets/S_{1}}{T_1}-multicategory on 1', or `free
\pr{\Sets}{\mr{free\ monoid}}-operad', on $\Sets/\nat$. If
$Y=(Y(n))_{n\elt\nat}$ is an object of $\Sets/\nat$, then a member $y$ of
$Y(n)$ can be drawn as
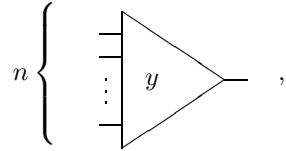
\[
n \left\{
\ctransistor{y}{}{}{}{}
\right.
,
\]
or for Baez-Dolan adherents, as
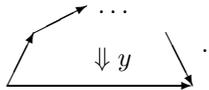
\[
\begin{opetope}
	&	&	&\cnr	&\ldots	&	&	&	\\
	&\cnr	&\ruEdge(2,1)&	&	&	&\cnr	&	\\
\ruEdge(1,2)&	&	&	&\Downarrow y&	&	&\rdEdge(1,2)\\
\cnr	&	&	&\rEdge	&	&	&	&\cnr	\\
\end{opetope}
\ .
\]
The monad $T_2$ sends $Y$ to the family of pictures obtained by sticking
together members of the $Y(n)$'s.

A description of opetopic sets can also be given, in a manner similar to that
for opetopes. Opetopic sets are central to the Baez-Dolan development
\cite{BD} of $n$-category theory; the explanation of opetopic sets most
convenient to us here is closer to that formulated in \cite{HMP}, as
interpreted to me by Martin Hyland from a conversation with John Power.

\subsection{The Grothendieck Construction, or Slicing}
\label{fds:slicing}

Given an ordinary category \scat{C} and a functor $\scat{C}\goby{h}\Sets$,
the Grothendieck construction produces a category $\scat{C}_{h}$ such that
\[
\ftrcat{\scat{C}_{h}}{\Sets}
\iso
\ftrcat{\scat{C}}{\Sets}/h.
\]
In general, given an \Cartpr-multicategory $A$ and an algebra
\[
\bktdvslob{X}{p}{S}^{\blob} \goby{h} \bktdvslob{X}{p}{S}
\]
for A (with \ubl\ as in Construction~\ref{constn:ind-mnd}), we may describe a
new \Cartpr-multicategory $A_h$ such that
\[
\Alg(A_{h}) = \Alg(A) / h.
\]
The graph of $A_h$ is
\begin{slopeydiag}
	&		&X_{\blob}	&		&	\\
	&\ldTo<{\phi_p}	&		&\rdTo>{h}	&	\\
X^{*}	&		&		&		&X	\\
\end{slopeydiag}
where $\phi_p$ is part of the pullback square defining $X_{\blob}$
(see~\ref{constn:ind-mnd}); identities and composition are given via the unit
and multiplication of the monad $\ubl$. (The natural map $A_h \go A$ is, in a
suitable sense, a discrete opfibration, and one can go on to show that
algebras for $A$ correspond exactly to discrete opfibrations over $A$.)

In \cite[section 2.5]{BD}, $A_h$ would be called a `slice operad'; slicing
plays an essential part in their theory.

It is perhaps worth noting that the slicing of multicategories corresponds to
the slicing of monads, in the following sense. Given a monad $T$ on a
category \scat{C} and an algebra $TC \goby{h} C$ for $T$, there is a monad
$T_h$ on $\scat{C}/C$ such that
\[
\Alg(T_h) = \Alg(T)/h,
\]
where $\Alg\blank$ denotes the category of algebras for a monad. Now suppose
we start with an \Cartpr-multicategory $A$ and an algebra $h$ for $A$, as
above. We get the monad \blbk\ on \Slice, and therefore a monad $\blbk_{h}$
on $\frac{\Slice}{X\go S} \iso \frac{\ess}{X}$. But we also get the
\Cartpr-multicategory $A_h$ on $X$, and therefore another monad on
$\ess/X$. The reader will not be surprised to learn that these two monads on
$\ess/X$ are the same.

\subsection{Structured Categories}	\label{fds:struc}

The observation from which this section takes off is that any strict monoidal
category has an underlying (ordinary) multicategory. (All monoidal categories
and maps between them will be strict in this section; one could consider
similar issues for lax versions, but this is not done here. For the time
being, `multicategory' means \pr{\Sets}{\mr{free\ monoid}}-multicategory.)
Explicitly, if \pr{\scat{C}}{\otimes} is a monoidal category, then the
underlying multicategory $A$ has the same object-set as \scat{C} and has
homsets defined by
\[
\mtihom{\Hom_{A}}{C_1}{C_n}{C} = 
\homset{\Hom_{\scat{C}}}{C_{1}\otimes\cdots\otimes C_{n}}{C}
\]
for objects \range{C_1}{C_n,C}. Composition and identities in $A$ are easily
defined.

There is a converse process: given any multicategory $A$ with objects $S$,
there is a `free' monoidal category \scat{C} on it.  Informally, an
object/arrow of \scat{C} is a sequence of objects/arrows of A.  Thus the
objects of \scat{C} are of form \bftuple{s_1}{s_n} ($s_i \elt S$), and a
typical arrow $\tuplebts{s_1,s_2,s_3,s_4,s_5} \go
\tuplebts{s'_1,s'_2,s'_3}$ is a sequence \tuplebts{a_1,a_2,a_3} of elements
of $A$ with domains and codomains as illustrated:
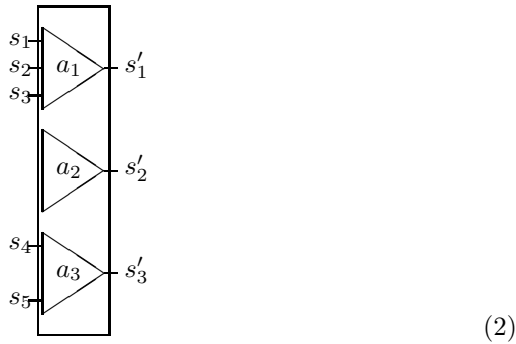
\begin{equation}	\label{diag:arrows-in-mon-cat} 
\freemoncatpic
\end{equation}
The tensor in \scat{C} is just juxtaposition.

For example, the terminal multicategory \fcat{1} has one object and, for each
$n\elt\nat$, one arrow of form
\[
n \left\{
\ctransistor{}{}{}{}{}
\right.
\ ;
\]
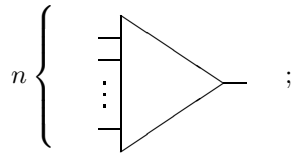
figure~\ref{diag:arrows-in-mon-cat} (above) indicates that the `free'
monoidal category on the multicategory \fcat{1} is $\Delta$, the
category of finite ordinals, with addition as $\otimes$.

The name `free' is justified: that is, there is an adjunction
\begin{diagram}
\mbox{Monoidal categories}	\\
\uTo\dashv\dTo			\\
\mbox{Multicategories}		\\
\end{diagram}
where the two functors are those described above. Moreover, this adjunction
is monadic. (But the forgetful functor does \emph{not} provide a full
embedding of Monoidal categories into Multicategories. It is faithful, but
not full: there is a multicategory map $\fcat{1} \go \Delta$ sending the
unique object of \fcat{1} to the object 1 of $\Delta$, and this map does not
preserve the monoidal structure.)

Naturally, we would like to generalize from $\Cartpr=\pr{\Sets}{\mr{free\
monoid}}$ to any cartesian \Cartpr. To do this, we need a notion of
`\Cartpr-structured category', which in the case \pr{\Sets}{\mr{free\
monoid}} just means monoidal category. One can view a monoidal category
either as a monoid in \fcat{Cat}, i.e.\ an algebra for the monoid monad on
\fcat{Cat}, or as a category object in \fcat{Monoids}. The latter view is
more convenient here: if $\ess^\stbk$ is the category of algebras for the
monad \stbk on \ess, then define an \emph{\Cartpr-structured category} to be
an \pr{\ess^\stbk}{\id}-multicategory, i.e.\ a category object in
$\ess^\stbk$.  (Alternatively, as in the motivating case, \ust\ can be
extended to give a monad on $\ess\fcat{-Cat}$, and an \Cartpr-structured
category defined as an algebra for this monad. It comes to the same thing.)

It is now possible to describe a monadic adjunction
\begin{diagram}
\Cartpr\mbox{-Structured categories}	\\
\uTo<{F}\dashv\dTo>{U}			\\
\Multicatof{\Cartpr}			\\
\end{diagram}
generalizing that above. The effect of the functors $F$ and $U$ on objects is
as outlined now. Given an \Cartpr-multicategory \graph{A}{S}{d}{c}, the
category $FA$ has
graph 
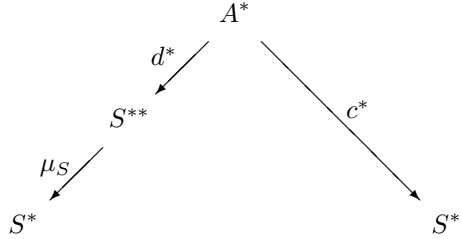
\begin{slopeydiag}
	&	&	&	&A^{*}	&	&	&	&	\\
	&	&	&\ldTo<{d^*}&	&\rdTo(4,4)>{c^*}&&	&	\\
	&	&S^{**}	&	&	&	&	&	&	\\
	&\ldTo<{\mu_S}&	&	&	&	&	&	&	\\
S^*	&	&	&	&	&	&	&	&S^*	\\
\end{slopeydiag}
and the monoidal structures $A^{**} \goby{\otimes} A^{*}$, $S^{**}
\goby{\otimes} S^{*}$ are components of $\mu$. Given an
\Cartpr-structured category \spaan{B}{R}{R}{q}{p}, with monoidal structure
$R^{*} \goby{\otimes} R$ and $B^{*} \goby{\otimes} B$, the graph \gph{A}{R}
of $UB$ is given by
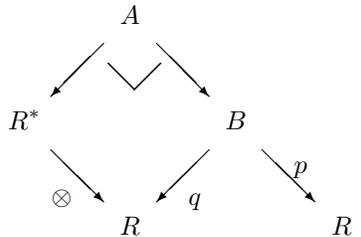
\[
\begin{slopeydiag}
	&		&A\Spbk	&	&	&	&	\\
	&\ldTo		&	&\rdTo	&	&	&	\\
R^{*}	&		&	&	&B	&	&	\\
	&\rdTo<{\otimes}&	&\ldTo>{q}&	&\rdTo>{p}&	\\
	&		&R	&	&	&	&R	\\
\end{slopeydiag}
\ .
\]

In fact, all of the above can be seen as a certain instance of the functorial
action of \fcat{Multicat}, as described in the next section.

\subsection{Functoriality of \Multicat}	\label{fds:multicat}

Any cartesian pair \Cartpr\ yields the category of \Cartpr-multicategories;
it would be reasonable to expect a map $\pr{\cat{R}}{\ubl} \go \Cartpr$ of
cartesian monads to yield a functor
\[
\Multicatof{\pr{\cat{R}}{\ubl}} \go \Multicatof{\Cartpr}
.
\]
As explained on page~\pageref{p:st-defns}, `map' might mean either monad
functor or monad opfunctor. Whichever meaning we take, we do get the kind of
functoriality desired, as long as the monad (op)functor is cartesian
(page~\pageref{p:cart-st-defns}). All this extends to 2-cells, so we have two
2-functors
\begin{eqnarray*}
\fcat{CartMnd} 	&\go 	&\Cat	\\
\fcat{CartMnd'}	&\go	&\Cat	\\
\end{eqnarray*}
agreeing on 0-cells.

We now sketch out how a cartesian monad (op)functor yields a functor between
multicategory categories. If $\pr{\cat{R}}{\ubl} \goby{\pr{P}{\phi}} \Cartpr$
is a cartesian monad functor, then 
\[
\Multicatof{\pr{\cat{R}}{\ubl}} \goby{\overline{P}} \Multicatof{\Cartpr}
\]
is defined by pullback: for an $\pr{\cat{R}}{\ubl}$-multicategory with graph
\spaan{B}{R^{\blob}}{R}{q}{p}, the graph of the multicategory $\overline{P}B$
is given by the diagram
\[
\begin{slopeydiag}
	&	&\overline{P}B\Spbk&	&	&	&	\\
	&\ldTo	&		&\rdTo	&	&	&	\\
(PR)^*	&	&		&	&PB	&	&	\\
	&\rdTo<{\phi_R}&	&\ldTo>{Pq}&	&\rdTo>{Pp}&	\\
	&	&P(R^{\blob})	&	&	&	&PR	\\
\end{slopeydiag}
\ \ .
\]
If $\pr{\cat{R}}{\ubl} \goby{\pr{P}{\phi}} \Cartpr$ is a cartesian monad
opfunctor, then $\overline{P}$ is defined by composition: $\overline{P}B$ has
graph
\[
\begin{slopeydiag}
	&	&	&	&PB	&	&	&	&	\\
	&	&	&\ldTo<{Pq}&	&\rdTo(4,4)>{Pp}&&	&	\\
	&	&P(R^{\blob})&	&	&	&	&	&	\\
	&\ldTo<{\phi_R}&&	&	&	&	&	&	\\
(PR)^*	&	&	&	&	&	&	&	&PR	\\
\end{slopeydiag}
\ \ .
\]

We have described two categories, \fcat{CartMnd} and \fcat{CartMnd'}, on
which \Multicat\ acts as a functor, but there is still another. Suppose we
have a diagram
\begin{diagram}[balance]
\pr{\cat{R}}{\ubl}				\\
	\uTo<{\mr{opfunctor\ }\pr{P}{\phi}}
	\dTo>{\mr{functor\ }\pr{Q}{\psi}}	\\
\Cartpr						\\
\end{diagram}
in which everything is cartesian, $P\dashv Q$ (as plain functors), and the
unit and counit of the adjunction commute suitably with $\phi$ and $\psi$.
Then there arises an adjunction
\[
\begin{diagram}[balance]
\Multicatof{\pr{\cat{R}}{\ubl}}	\\
	\uTo<{\overline{P}}
	\dashv	
	\dTo>{\overline{Q}}	\\
\Multicatof{\Cartpr}		\\
\end{diagram}
\ \ ,
\]
defined in an evident way (with $\overline{P}$ and $\overline{Q}$ as
above). In particular, we can apply this to
\begin{diagram}[balance]
\pr{\ess^{\stbk}}{\id}		\\
	\uTo<{\pr{F}{\phi}}
	\dTo>{\pr{U}{\psi}}	\\
\Cartpr				\\
\end{diagram}
for any cartesian \Cartpr, where $\ess^{\stbk}$ is the category of algebras
for the monad \ust\ on \ess, $F$ and $U$ are the free algebra and forgetful
functors, and $\phi$ and $\psi$ are certain canonical natural
transformations. This gives the adjunction
\begin{diagram}[balance]
\Cartpr\mbox{-Structured categories}	\\
\uTo\dashv\dTo				\\
\Multicatof{\Cartpr}			\\
\end{diagram}
of~\ref{fds:struc}.

\pagebreak
\subsection{The Endomorphism Multicategory}	\label{fds:endo}

Any set $X$ gives rise to an operad $E=\End(X)$ (for the free monoid monad on
\Sets); it is defined by
\[
E(n)=\Sets\pr{X^n}{X},
\]
with evident units and composition functions. (Recall that for us, an operad
is a multicategory with just one object.) Given any operad $A$, one may
define an algebra for $A$ to be a set $X$ together with an operad map $A \go
\End(X)$, and, of course, this is equivalent to the definition of algebra
given above. Many theories of operads, e.g. \cite{Bat}, define `algebra' in
this fashion, so we indicate here how it fits into the general theory.

Suppose we have an \Cartpr-multicategory \gph{A}{S}, and that the category
$\frac{\ess}{S^{*}\times S}$ of \Cartpr-graphs on $S$ has exponentials. This
occurs, for instance, if \ess\ is a topos. Let \slob{X}{p}{S} be an object of
\Slice, and put
\[
E=	\left[
	\vslob{X^{*}\times S}{p^{*}\times 1}{S^{*}\times S}
	,
	\vslob{S^{*}\times X}{1\times p}{S^{*}\times S}
	\right]
\]
where \ftrcat{\emptybk}{\emptybk} indicates exponential. Then $E$ carries a
natural \Cartpr-multicategory structure, and algebra structures on $X$
correspond to multicategory maps $A \go E$.

\subsection{Characterization of Multicategories by Monads}
\label{fds:characterization}

A traditional May-style operad induces a monad on \Sets, whose algebras are
the algebras of the operad. In our general setting, an \Cartpr-multicategory
on $S$ induces a monad on \Slice. In both cases, one may ask precisely which
monads arise from multicategories/operads. Here we answer the question by
giving a complete description of multicategories in terms of monads
(Lemma~\ref{lemma:mnd-data}). 

\begin{llemma}	\label{lemma:multi-to-mnd}
Let \Cartpr\ be cartesian, \gph{A}{S} an \Cartpr-multicategory, and \ubl\ the
induced monad on \Slice. Then the forgetful functor $\Slice\goby{U}\ess$
naturally carries the structure of a monad opfunctor, this opfunctor is
cartesian, and \ubl\ is a cartesian monad.
\end{llemma}

(For the definition of cartesian monad opfunctor, see
page~\pageref{p:st-defns} ff.)

\pf\ The data required is a natural transformation
\[
\begin{ntdiag}
\Slice	&	&\rTo^{\blob}	&	&\Slice		\\
	&	&		&\	&		\\
\dTo<{U}&	&\ldTo>{\phi}	&	&\dTo>{U}	\\
	&\	&		&	&		\\
\ess	&	&\rTo_{*}	&	&\ess		\\
\end{ntdiag}
\ .
\]
That is, for each object \slob{X}{p}{S} of \Slice, we need a map $X_{\blob}
\goby{\phi_p} X^*$; this map is part of the pullback square defining
$X_{\blob}$ in~\ref{constn:ind-mnd}. The rest of the proof is easy checking.
\done 

In fact, this monad data arising from $A$ characterizes completely
\Cartpr-multicategories on $S$:

\begin{llemma}	\label{lemma:mnd-data}
Let \Cartpr\ be cartesian and $S\elt\ess$. Then an \Cartpr-multicategory on
$S$ is the same thing as a cartesian monad on \Slice\ together with the
structure of a cartesian monad opfunctor on $\Slice\goby{U}\ess$.
\end{llemma}

\pf\ Lemma~\ref{lemma:multi-to-mnd} shows how an \Cartpr-multicategory $A$
yields the monad data $\Icartpr \goby{\pr{U}{\phi}} \Cartpr$. It is easy to
see that $\bktdslob{S}{1}{S}^{\blob} = \slob{A}{c}{S}$ and that
$\phi_{\slob{S}{1}{S}}$ is $A \goby{d} S^*$, so from this
monad data we can recover the graph structure of the
multicategory. Similarly, \ids\ and \comp\ are recovered as the unit and
multiplication of \ubl\ at \slob{S}{1}{S}. This tells us how to pass from
monad data to a multicategory. \done

As an application of this result, consider the strongly regular algebraic
theories~(\ref{eg:cart-mnds}(\ref{eg:sr-theories})).  If \ubl\ is the monad
on \Sets\ from some strongly regular theory and \ust\ the free monoid monad,
then any strongly regular presentation of the theory gives rise to a
cartesian natural transformation $\ubl\go\ust$, which commutes with the monad
structure. We therefore have a cartesian monad \ubl\ on $\Sets=\Sets/1$, and
the structure of a cartesian monad opfunctor on the forgetful functor
$\Sets/1\go\Sets$. By Lemma~\ref{lemma:mnd-data}, \ubl\ arises from a
\pr{\Sets}{\mr{free\ monoid}}-multicategory on 1. Thus, given a strongly
regular theory, there is an (ordinary) operad whose algebras are the same as
those of the theory.

Lemma~\ref{lemma:multi-to-mnd} says in particular that the monad \ubl\ on
\Slice\ arising from an \Cartpr-multicategory $A$ on $S$ is cartesian; one
may therefore ask what an \Icartpr-multicategory is. The answer is simple:
\[
\Multicatof{\Icartpr}	\iso	\Multicatof{\Cartpr}/A.
\]

\subsection{A Bicategorical Description}	\label{fds:bicat}

In this section I will give an alternative definition of
\Cartpr-multicategory, using weak 2-monads on bicategories. The significance
of this description eludes me, and a notable omission is a definition in this
framework of an algebra for a multicategory.

Let \Cartpr\ be cartesian. Then there is a kind of weak 2-monad \ulid\
induced on the bicategory $\fcat{Spans}_{\ess}$ of spans in \ess, as follows.
The `functor' part is illustrated by the picture
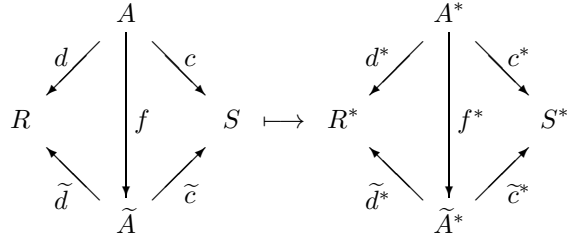
\[
\begin{slopeydiag}
 &			&A		&			&\\
 &\ldTo<{d}		&		&\rdTo>{c}		&\\
R&			&\dTo>{f}	&			&S\\
 &\luTo<{\twid{d}}	&		&\ruTo>{\twid{c}}	&\\
 &			&\twid{A}	&			&\\
\end{slopeydiag}
\ \ \goesto\ \  
\begin{slopeydiag}
   &			&A^*		&			&\\
   &\ldTo<{d^*}		&		&\rdTo>{c^*}		&\\
R^*&			&\dTo>{f^*}	&			&S^*\\
   &\luTo<{\twid{d}^*}	&		&\ruTo>{\twid{c}^*}	&\\
   &			&\twid{A}^*	&			&\\
\end{slopeydiag}
\]
for 2-cells. This `functor' preserves 1-cell composition up to isomorphism;
that is, it is a homomorphism of bicategories. The rest of the `monad'
structure is described in a similar manner, and it turns out that what we
have is:
\begin{itemize}
\item 	a homomorphism $\fcat{Spans}_{\ess} \goby{\lidbk} \fcat{Spans}_{\ess}$
\item 	strong transformations $1 \goby{\eta} \lidbk \ogby{\mu}
	\blank^{\lid\lid}$ 
\end{itemize}
such that the monad axioms are satisfied up to invertible modification. (The
meaning of these technical terms is that all the 1-cell diagrams stating
naturality, associativity, etc., hold up to isomorphism.) Such a structure on
a bicategory will just be called a \emph{``monad''}, in inverted commas.

Given a ``monad'' \ulid\ on a bicategory \cat{B}, define an
\emph{``algebra''} for \ulid\ to be a 0-cell $S$ together
with a 1-cell $S^{\lid}\goby{h}S$ and 2-cells
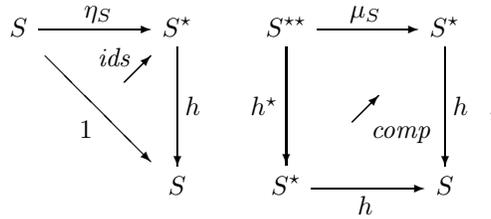
\[
\begin{ntdiag}
S	&	&\rTo^{\eta_S}	&		&S^{\lid}	\\
	&\rdTo(4,4)<{1}&	&\ruTo<{\ids}	&		\\
	&	&		&		&\dTo>{h}	\\
	&	&		&		&		\\
	&	&		&		&S		\\
\end{ntdiag}
\ \ \ \ \ 
\begin{ntdiag}
S^{\lid\lid}	&	&\rTo^{\mu_S}	&	&S^{\lid}	\\
		&	&		&\	&		\\
\dTo<{h^{\lid}}	&	&\ruTo>{\comp}	&	&\dTo>{h}	\\
		&\	&		&	&		\\
S^{\lid}	&	&\rTo_{h}	&	&S		\\
\end{ntdiag}
\ \ ,
\]
such that the 2-cells satisfy equations looking like associativity and
identity laws.

An \Cartpr-multicategory is then the same thing as an ``algebra'' for the
``monad'' \ulid\ on $\fcat{Spans}_{\ess}$: the 1-cell $S^{\lid} \goby{h} S$
is a diagram \graph{A}{S}{d}{c} in \ess, and the 2-cells \ids\ and \comp\
have the roles suggested by their names.


\clearpage



\end{document}